\documentclass[reqno]{amsart}
\makeatletter
\@namedef{subjclassname@2020}{%
  \textup{2020} Mathematics Subject Classification}
\makeatother
\usepackage{amsmath,amssymb}
\usepackage{amscd}
\usepackage[none]{hyphenat}
\usepackage{soul}
\usepackage[utf8]{inputenc} 
\usepackage[T1]{fontenc}    
\usepackage{cancel}
\usepackage{caption}
\usepackage{url}            
\usepackage{booktabs}       
\usepackage{amsfonts}       
\usepackage{nicefrac}       
\usepackage{microtype}      
\usepackage{color}
\usepackage{lipsum}
\usepackage[T1]{fontenc}
\usepackage{float}
\usepackage{calrsfs}
\usepackage{enumerate}
\usepackage{tikz-cd}
\usepackage{tikz}
\tikzcdset{scale cd/.style={every label/.append style={scale=#1}, cells={nodes={scale=#1}}}}
\usepackage{graphicx}
\usepackage{old-arrows} 
\usetikzlibrary{arrows}
\usepackage{comment}
\usepackage{blindtext}
\usepackage[mathscr]{euscript}
\DeclareMathAlphabet{\pazocal}{OMS}{zplm}{m}{n}

\usepackage{eqnarray,amsmath}
\usepackage{marginnote}
\usepackage{mathtools}
\usepackage{amsbsy}
\usepackage[hidelinks]{hyperref}
\usepackage{animate} 

\input xygraph
\input xy
\xyoption{all}

\usepackage{makecell,multirow}
\usepackage{rotating}
\newcolumntype{C}[1]{>{\centering\let\newline\\\arraybackslash\hspace{0pt}}m{#1}}

\usepackage[many]{tcolorbox}
\newtheorem{lemma}{Lemma}[section]

\newtheorem{theorem}[lemma]{Theorem}
\newtheorem{proposition}[lemma]{Proposition}
\newtheorem*{proposition*}{Proposition}
\newtheorem*{corollary*}{Corollary}
\newtheorem*{theorem*}{Theorem}
\newtheorem*{lemma*}{Lemma}

\theoremstyle{definition}
\newtheorem{definition}[lemma]{Definition}

\newtheorem{example}[lemma]{Example}

\newtheorem{remark}[lemma]{Remark}

\author[D. Mart\'\i n]{Dolores Mart\'\i n Barquero}
\address{D. Mart\'\i n Barquero: Departamento de Matem\'atica Aplicada, Escuela de Ingenier\'\i as Industriales, Universidad de M\'alaga, Campus de Teatinos s/n. 29071 M\'alaga.   Spain.}
\email{dmartin@uma.es}

\author[C. Mart\'\i n]{C\'andido Mart\'\i n Gonz\'alez}
\address{C. Mart\'\i n Gonz\'alez: Departamento de \'Algebra Geometr\'{\i}a y Topolog\'{\i}a, Fa\-cultad de Ciencias, Universidad de M\'alaga, Campus de Teatinos s/n. 29071 M\'alaga.   Spain.}
\email{candido\_m@uma.es}

\author[I. Ruiz]{Iv\'an Ruiz Campos}
\address{I. Ruiz Campos:  Departamento de \'Algebra Geometr\'{\i}a y Topolog\'{\i}a, Fa\-cultad de Ciencias, Universidad de M\'alaga, Campus de Teatinos s/n. 29071 M\'alaga. Spain.}
\email{ivaruicam@uma.es}

\newcommand{\RR}{\mathbb{R}}

\newcommand{\MM}{\mathcal{M}}
\newcommand{\F}{\mathcal{F}}
\newcommand{\G}{\mathcal{G}}

\newcommand{\C}{\mathcal{C}}

\newcommand{\Vect}{\textbf{Vect}}

\def\Dg{\hbox{\bf DiGr}}

\DeclareMathOperator{\PP}{Path}

\DeclareMathOperator{\Span}{span}
\DeclareMathOperator{\Reg}{Reg}

\DeclareMathOperator{\Id}{Id}

\DeclareMathOperator{\h}{h_{alg}}
\DeclareMathOperator{\Falg}{\mathbf{Falg}_{\mathit{K}}}
\DeclareMathOperator{\Alg}{\mathbf{Alg_{K}}}
\DeclareMathOperator{\DG}{\mathbf{DiGr}}

\DeclareMathOperator{\CKD}{\mathbf{CKGr}}
\DeclareMathOperator{\falg}{\mathbf{falg}_{\mathit{K}}}

\DeclareMathOperator{\tr}{tr}

\DeclareMathOperator{\spec}{spec}

\def\s{\sigma}

\def\h{\mathop{\hbox{\rm h}_{\hbox{\tiny\rm alg}}}}

\def\remove#1{}
\def\F{\mathcal{F}}
\def\G{\mathcal{G}}

\title{An extension of the notion of algebraic entropy}

\begin{document}
\begin{abstract}
  The definition of algebraic entropy presents a difficulty in the general setting of filtered algebras, which makes it necessary to distinguish between filtered algebras and finitely filtered algebras. While algebraic entropy is well defined in the finitely filtered case, extending this notion to arbitrary filtered algebras is not immediate. To address this issue, we realize a filtered algebra as the limit of a strict direct system of finitely filtered algebras and define its entropy as the supremum of the entropies of the algebras in the system. This extension is consistent with the original definition, when the algebra is finitely filtered, the extended entropy agrees with its algebraic entropy. This new setting allows us to extend the definition of the algebraic entropy to path and Leavitt path algebras of infinite graphs.
\end{abstract}
\maketitle

\section{Introduction}
In the study of infinite dimensional algebras, the notion of dimension does not contribute much information. However, instead of the dimension, one can categorize an algebra by the growth of the dimension of finitely generated subalgebras, if the algebra is filtered or graded. It is somewhat unsatisfactory that a quantification of the exponential growth is not given via this concept. For this purpose, in the case of groups and filtered algebras, the notion of algebraic entropy is introduced, see e.g.~ \cite{dikranjan2012connection,newman2000entropy}. The concept of entropy itself was sketched for the first time by Adler, Konheim and McAndrew for compact topological spaces \cite{Adler}. They define the topological entropy via covers and apply it to dynamical systems and endomorphisms of abelian groups. A generalisation of this notion is given by D. Gonçalves, D. Royer and F. A. Tasca in \cite{entropydaniel} which is used in infinite alphabet shift spaces. This algebraic entropy, was later studied by Weiss \cite{Weiss} using Pontryagin duality.
A different definition was given by Peters \cite{peters1979} for automorphisms of abelian groups.  It coincides with the one by Weiss on endomorphisms of torsion abelian groups. The algebraic entropy was generalized in group theory in e.g.~\cite{dikranjan2016entropy,bruno2017some}. There are also contributions to the definition of algebraic entropy for semigroup and monoids \cite{entropychinos}.
For graded algebras the entropy was studied in \cite{newman2000entropy}. The growth of the algebra is in this setting based on the dimensions of certain factor spaces. In statistical mechanics, the entropy consists of counting the number of possible states a system can be in. The algebraic entropy can be considered similarly. The microstates that Boltzmann \cite{Boltzmann} used in his definition of the entropy are replaced by the span of finite dimensional subspaces generating the infinite dimensional algebra. The growth is then defined as the dimension of consecutive factor spaces. In the present article we extend this concept to general path algebras and Leavitt path algebras. The Leavitt path algebra of a directed graph in the past decades received substantial interests from algebraists after its introduction in two separate, but simultaneously created, papers by Abrams and Aranda \cite{AbramsPino} and by Ara, Moreno and Pardo \cite{AraMorenoPardo}. Studies look for description of the algebraic structure of these algebras as well as classifications via the geometry of the associated graphs. In \cite{moreno2018graph} the Gelfand-Kirillov dimension (GK-dimension) of the path algebras were investigated whereas in \cite{Zelmanov12} for the Leavitt path algebra. By generalizing the entropy of graded algebras given by \cite{newman2000entropy} the concept of algebraic entropy for filtered algebras is defined in \cite{entropy1} and applied to path and Leavitt path algebras of finite graphs. Latter, the conjecture of the two values of entropy coinciding for their corresponding standard filtration is proven in \cite{entropy2}. However, it was still missing a definition of entropy that can be applied to infinite graphs. Following \cite{AraMorenoPardo,Goodearl}, it is known that every infinite graph is the direct limit of its finite subgraphs. Furthermore, when the graph is row-finite and the subgraphs are complete, the path algebra and the Leavitt path algebra of the graph are the direct limits of the respective path algebras and Leavitt path algebras associated with its complete subgraphs. 
These give us hope of defining a concept of entropy for limits of directed systems of filtered algebras.
In this sense, we obtain a generalization of the algebraic entropy given in \cite{entropy1,entropy2}.

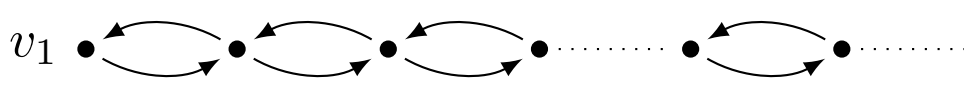
\begin{figure}[h]
    \centering
    \begin{tikzpicture}[shorten <=2pt,shorten >=2pt,>=latex, node distance={20mm}, main/.style = {draw, fill, circle, inner sep = 1pt}]

    \node[main,label = left:$v_1$] (1) at (0,0) {};
    \node[main] (2) at (1,0) {};
    \node[main] (3) at (2,0) {};
    \node[main] (4) at (3,0) {};
    \node[main] (5) at (4,0) {};
    \node[main] (6) at (5,0) {};
    \node (7) at (6,0) {};

    \draw[->] (1) to [bend right] (2);
    \draw[->] (2) to [bend right] (3);
    \draw[->] (2) to [bend right] (1);
    \draw[->] (3) to [bend right] (2);
    \draw[->] (3) to [bend right] (4);
    \draw[->] (4) to [bend right] (3);
    \draw[dotted] (4) to (5);
    
    \draw[->] (5) to [bend right] (6);
    \draw[->] (6) to [bend right] (5);

    \draw[dotted] (6) to (7); 
    \end{tikzpicture}
    \caption[Infinity chain]{Infinity chain.}
    \label{fig_chorizo}
\end{figure}

To start this work, we begin with a motivating example.
Following \cite{entropy1,entropy2}, denote the graph $Ch_\infty$ given by vertices $Ch_\infty^0 = \{v_n\}_{n=1}^\infty$ and edges $Ch_\infty^1 = \{e_n,f_n\}_{n=1}^\infty$ with $s(e_n)= r(f_n) = v_n $ and $r(e_n) = s(f_n) = v_{n+1}$. We denote it as the \emph{infinity chain} and it can be represented as in Figure \ref{fig_chorizo}.
Consider the path algebra $KCh_\infty$, if we think of its standard filtration $\F=\{V_n\}_{n=1}^\infty $, given by $ V_n = \Span_K \{\lambda \in \PP(E) \colon | \lambda|\leq n\}$, we can observe that $\dim(V_n/V_{n-1}) = \left |\{\lambda \in \PP(E) \colon |\lambda| = n\} \right | = +\infty$ which means that 
\begin{equation*}
    \h(KCh_\infty) = \limsup_{n\to \infty} \frac{\dim \left(V_n/V_{n-1} \right )}{n} = +\infty. 
\end{equation*}

\begin{figure}[h]
    \centering

    \begin{tikzpicture}[
        shorten <=2pt, shorten >=2pt, >=latex,
        node distance={15mm},
        main/.style={draw, fill, circle, inner sep=1pt}
    ]
        \node[main, label=left:$v_1$] (1) at (0,0) {};
        \node[main] (2) at (1,0) {};
        \node[main] (3) at (2,0) {};

        \draw[->] (1) to[bend right] (2);
        \draw[->] (2) to[bend right] (3);
        \draw[->] (2) to[bend right] (1);
        \draw[->] (3) to[bend right] (2);
    \end{tikzpicture}
    
        $\vdots$
\vspace{2mm}

    \begin{tikzpicture}[
        shorten <=2pt, shorten >=2pt, >=latex,
        node distance={15mm},
        main/.style={draw, fill, circle, inner sep=1pt}
    ]
        \node[main, label=left:$v_1$] (1) at (0,0) {};
        \node[main] (2) at (1,0) {};
        \node[main] (3) at (2,0) {};
        \node[main] (4) at (3,0) {};

        \draw[->] (1) to[bend right] (2);
        \draw[->] (2) to[bend right] (1);
        \draw[->] (2) to[bend right] (3);
        \draw[->] (3) to[bend right] (2);
        \draw[->] (3) to[bend right] (4);
        \draw[->] (4) to[bend right] (3);
    \end{tikzpicture}

        $\vdots$
\vspace{2mm}

    \begin{tikzpicture}[
        shorten <=2pt, shorten >=2pt, >=latex,
        node distance={15mm},
        main/.style={draw, fill, circle, inner sep=1pt}
    ]
        \node[main, label=left:$v_1$] (1) at (0,0) {};
        \node[main] (2) at (1,0) {};
        \node at (2,0) {$\cdots$};
        \node[main] (n-1) at (3,0) {};
        \node[main, label=right:$v_n$] (n) at (4,0) {};

        \draw[->] (1) to[bend right] (2);
        \draw[->] (2) to[bend right] (1);

        \draw[->] (n-1) to[bend right] (n);
        \draw[->] (n) to[bend right] (n-1);
    \end{tikzpicture}
   
        $\vdots$
\vspace{1mm}
    \caption{Representations of the chains $Ch_2$, $Ch_3$ and
    the general graph $Ch_n$.}
    \label{fig_chor}
\end{figure}
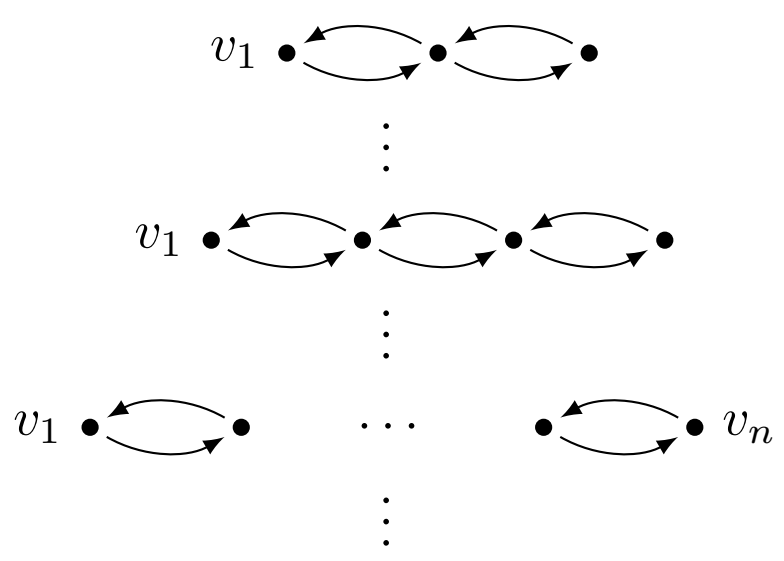

 On the other hand, consider the family of subgraphs $\{Ch_n\}_{n = 2}^\infty$ given by Figure \ref{fig_chor}. The reader can observe that their corresponding adjacency matrices are
\begin{equation*}
    A_n = \tiny \left (\begin{array}{ccccccc}
        0 & 1 & 0 & \cdots & 0& 0 \\
        1 & 0 & 1 & \cdots & 0& 0 \\
        0 & 1 & 0 & \cdots & 0&0 \\
        \vdots & \vdots &\vdots & \ddots & \vdots & \vdots \\
        0 & 0 & 0 & \cdots & 0& 1 \\
        0 & 0 & 0 & \cdots & 1 & 0
    \end{array} \right ) \in \MM_{n+1}(\RR).
\end{equation*}

Using the software \emph{Mathematica} we have computed numerically the values $\h(Kh_n)$ for several $n$ and represented them in Table \ref{tab_entropy_chorizo}.
\begin{table}[h]
\centering
\begin{tabular}{|c|c|}
\hline
$\mathbf{n}$ & $\mathbf{\h(KCh_n)}$ \\ \hline
$5$          & $0.549306$          \\ \hline
$15$         & $0.673745$          \\ \hline
$105$        & $0.692708$          \\ \hline
$203$        & $0.693029$          \\ \hline
$275$        & $0.693082$          \\ \hline
\end{tabular}
\caption{Algebraic entropy of the family of graphs $Ch_n$ for different values of $n$.}
\label{tab_entropy_chorizo}
\end{table}
The reader can observe that there is a certain tendency to the value $\log(2) \simeq 0.693147$ which is quite opposite to the value of entropy obtained with the definition given by \cite{entropy1}.

As it is natural, if we want to extend the definition of entropy of filtered algebras to algebras of infinite graphs, it has to be compatible with the one given for finite graphs \cite{entropy1}. For this purpose, we add a new condition to the definition of filtration and use the concept of entropy given in \cite{entropy1} for these ones. The work is organized as follows. In Section \ref{sec_preliminaries}, we set the basis of this theory. We give the notion of filtration and introduce the concept of finite filtration in Definition \ref{def_entropy_filtracion}. Furthermore, we distinguish between $\Falg$, the category of filtered algebras, and its subcategory $\falg$ of finitely filtered algebras, as defined in Definition \ref{def_entropy_filteredalgebra_cat}. In addition, we introduce the well notion of directed system and give the construction of it for the category of algebras and vector spaces. Since we are working with graphs, we make echo of the results in \cite{AraMorenoPardo,Goodearl} to study the limit of directed graphs, path and Leavitt path algebras. 

Section \ref{sec_limit_filtration} is dedicated to study direct limits in $\Falg$ and $\falg$. We prove in Proposition \ref{prop_limit_filteredalgebras} that the first one is closed under direct limits. However, the category of finitely filtered algebras is not, as we prove in Proposition \ref{prop_limitefinito_nofinito}.

Later, in Section \ref{sec_entropy_directsystems}, we focus our attention on the notion of a strict direct system (Definition \ref{def_entropy_strict}), for which we construct a new notion of entropy for algebras, that are limits of finitely generated algebras, in Definition \ref{def_entropy_limit}. Besides, we compute the value of the entropy of the path algebra that introduces and motivates this work. As it happens in \cite{entropy1} the value of the entropy depends on the chosen directed system. This is a direct consequence of Proposition \ref{prop_entropy_limit_dependencia}. We conclude this section with Subsection \ref{sub_compatibility} where we compare the new definition of entropy with the one given in \cite{entropy1} and verify that both are compatible in the sense of Theorem \ref{th_entropy_limit_filtered}.

In Section \ref{sec_finitely_generated}, we study the entropy for filtered algebras that are the direct limit of finitely generated finitely filtered algebras. In this context, the value of the entropy does not depend on the strict directed system used rather than just the original filtered algebra. We prove this result in Theorem \ref{th_independencia_entropía}. 

We end this work with Section \ref{sec_further_examples}. Since path and Leavitt path algebras are direct limits of finitely generated path and Leavitt path algebras, we prove that the compatibility of the values of entropies for these algebras (as in \cite{entropy2}) remains true for arbitrary row-finite graphs (Proposition \ref{prop_igualdad}). Finally, we give an examples for which we compute the entropy. 

\section{Preliminaries}\label{sec_preliminaries}
A \emph{directed graph} $E=(E^0,E^1,r,s)$ consists of two sets $E^0$ and $E^1$  and functions $r,s\colon E^1 \to E^0$. The elements of $E^0$ are called \emph{vertices} and the elements of $E^1$ \emph{edges}. We will mention that a graph $E$ is finite (resp. countable)  if $\vert E^0 \cup E^1 \vert < \infty$ where $\vert \cdot \vert$ denotes the cardinal of the corresponding set (resp. $\vert E^0 \cup E^1 \vert \leq  \aleph_0$). If a vertex $v$ emits no edges, that is, if $s^{-1}(v)$ is empty, then $v$is called a \emph{sink}. 
A vertex $v$ is called a \emph{regular vertex} if $s^{-1}(v)$ is a finite non-empty set. The set of regular vertices is denoted by ${\rm Reg}(E)$. A \emph{path of length $n$} in $E$ is a sequence $e_1 e_2 \ldots e_n$ of edges in $E$ such that $r(e_i)=s(e_{i+1})$ for $i \in \{1,2,\ldots,n-1\}$. If $\xi$ is a path of length $n$, then we write $|\xi|=n$.  We consider vertices in $E^0$ as paths of length zero. The set of all finite paths of length $n$ is denoted by $E^n$ and we let ${\rm Path}(E)=\cup_{n=0}^{\infty} E^n$. For instance,
for a field $K$ consider the \emph{path algebra} $KE$ as the free $K$-algebra generated by 
the sets $E^0 \cup E^1$ with relations:
\begin{itemize}
    \item[(V)] $vw=\delta_{v,w}v$ for $v,w \in E^0$,
    \item[(E1)] $s(e)e=er(e)=e$ for $e \in E^1$,
\end{itemize}
that is, $KE=\left \{ \sum_{i=1}^n k_i \lambda_i : k_i \in K, \lambda_i \in \text{Path}(E)\right \}$. 
The \emph{extended graph} of $E$ is defined as the new graph $\widehat{E}=(E^0, E^1\cup (E^1)^*, r',s'  )$, where $(E^1)^*= \{ e^*:e\in E^1 \}$ and the maps $r'$ and $s'$ are defined as
$r'|_{E^1}=r$, $s'|_{E^1}=s$, $r'(e^*)=s(e)$, and $s'(e^*)=r(e)$ for all $e\in E^1$. In other words, each $e^*\in (E^1)^*$ has orientation the reverse of that of its counterpart $e\in E^1$. The elements $(E^1)^*$ are called $ghost$ $edges$.

For a graph $E$ and a filed $K$, the  \textit{Leavitt path algebra} of $E$, denoted by $L_R(E)$, is the path algebra over the extended graph $\hat{E}$ with additional relations
\begin{enumerate}
\item[\textnormal{(CK1)}] $e^*e'=\delta_{e,e'}r(e)$ for all $e, e'\in E^1$;

\item[\textnormal{(CK2)}]$\sum_{ \{e \in E^1 \colon s(e)=v    \}  } e e^*=v$ for every $v\in \text{Reg}(E)$.
\end{enumerate}

\begin{figure}[h]
    \centering
    \begin{tikzpicture}[shorten <=2pt,shorten >=2pt,>=latex, node distance={15mm}, main/.style = {draw, fill, circle, inner sep = 1pt}]

    \node[main,label = left:$v_1$] (1) at (0,0) {};
    \node[main] (2) at (1,0) {};
    \draw[->] (1) to [bend right] (2);
    \end{tikzpicture}

    \begin{tikzpicture}[shorten <=2pt,shorten >=2pt,>=latex, node distance={15mm}, main/.style = {draw, fill, circle, inner sep = 1pt}]

    \node[main,label = left:$v_1$] (1) at (0,0) {};
    \node[main] (2) at (1,0) {};
    \node[main] (3) at (2,0) {};
    \draw[->] (1) to [bend right] (2);
    \draw[->] (2) to [bend right] (3);
    \draw[->] (2) to [bend right] (1);
    \end{tikzpicture}

    \begin{tikzpicture}[shorten <=2pt,shorten >=2pt,>=latex, node distance={15mm}, main/.style = {draw, fill, circle, inner sep = 1pt}]

    \node[main,label = left:$v_1$] (1) at (0,0) {};
    \node[main] (2) at (1,0) {};
    \node[main] (3) at (2,0) {};
    \node[main] (4) at (3,0) {};

    \draw[->] (1) to [bend right] (2);
    \draw[->] (2) to [bend right] (3);
    \draw[->] (2) to [bend right] (1);
    \draw[->] (3) to [bend right] (2);
    \draw[->] (3) to [bend right] (4);
    \end{tikzpicture}
    \caption[Chain of three links]{Some complete subgraphs of $Ch_\infty$.}
    \label{fig_chorr_completo}
\end{figure}
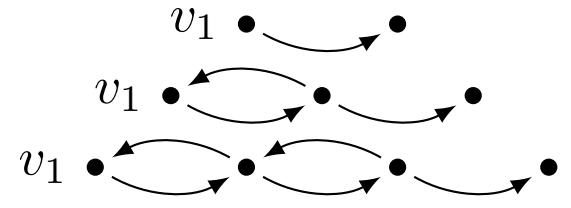

Considering a directed graph $E = (E^0,E^1,r,s)$ a \emph{subgraph} $F$ of $E$ is the directed graph $F = (F^0,F^1,r_F,s_F)$, such that $F^0\subseteq E^0$, $F^1\subseteq E^1$, and $r_F = \left. r \right |_{F^1}$ and $s_F = \left. s \right |_{F^1}$. We say that the subgraph $F$ is \emph{complete} if $s_F^{-1}(v) = s^{-1}(v)$ for all $v \in \Reg(F) \subseteq E^0$. In Figure \ref{fig_chorr_completo} we have some examples of complete subgraphs of $Ch_\infty$ 

\begin{definition}\label{def_entropy_filtracion}
    A $K$-algebra $A$ is said to be \emph{filtered}\index{algebra!filtered} if it is endowed with a collection of subspaces $\mathcal{F}=\{V_n\}_{n=0}^\infty$ such that
\begin{enumerate}
\item  $0 = V_0 \subseteq V_1\subseteq \cdots\subseteq V_n\subseteq V_{n+1}\subseteq \cdots \subseteq A$,  
\item $A=\bigcup_{n\ge 0}V_n$,
\item $V_nV_m\subseteq V_{n+m}$.
\end{enumerate}
The collection of vector spaces $\F$ is called a \emph{filtration}\index{filtration} of $A$. We say that the filtration is \emph{finite} if it also verifies 
\begin{enumerate}
    \item[(4)] $\dim(V_n/V_{n-1}) <+\infty$ for every $n\geq 1$. 
\end{enumerate}
\end{definition}

\begin{definition}\label{def_entropy_filteredalgebra_cat}
     We consider the \emph{category of filtered (resp. finitely filtered)}\index{category!of filtered algebras} $K$-algebras, denoted by $\Falg$ (resp. $\falg$). Its objects are the couples $(A,\F)$ where $\F$ is a filtration (resp. finite filtration) in $A$ and its morphisms $(A,\F)$ to $(B,\G)$ are the $K$-algebra morphisms $f\colon A\to B$ such that $f(V_n)\subseteq W_n$, for all $V_n \in \F$ and $W_n \in \G$. 
\end{definition}

Hence, following \cite{entropy1} the \emph{algebraic entropy of a finitely filtered algebra} $(A,\F)$ is $$\h(A,\F):=
\begin{cases}
0 &\text{ if } $A$ \text{ is finite dimensional,} \\
\displaystyle \limsup_{n\to\infty}\frac{
\log\dim(V_n/V_{n-1})}{n} & \text{ otherwise.}
\end{cases}$$
\begin{definition}
    If $E$ is a finite directed graph, then we define the \textit{standard filtration} of $KE$ (resp. $L_K(E))$ as the family of vector spaces $\F = \{V_n\}_{n\geq 1}$ (resp. $\G = \{W_n\}_{n \geq 1}$) given by
    \begin{equation*}
        V_n = \Span_K\{\lambda \in \PP(E) \colon |\lambda| \leq n\}
    \end{equation*}
   \begin{equation*}
        W_n = \Span_K\{\lambda\mu^* \in \PP(E)\PP(E)^* \colon |\lambda| + |\mu| \leq n\}.
   \end{equation*} 
\end{definition}

 A \emph{directed set} is a nonempty set $I$ together with a preorder $\le$ 
such that for every $i, j \in I$, there exists a $k \in I$ with 
$i \le k$ and $j \le k$. There is a canonical way to view a directed set 
$I$ as a small category. From this perspective, \emph{a directed system} in a 
category $\mathcal{C}$ is simply a functor $S \colon I \to \mathcal{C}$ (also called a \emph{diagram}), 
and its \emph{direct limit} $\varinjlim S$ (if it exists) is defined as the initial object in the category of cocones (Defined in \ref{def_direct_systems})  under $S$. However, we prefer a more down-to-earth 
approach, which we outline below.

\begin{definition}\label{def_direct_systems}
    Consider a category $\C$ and $I$ a directed set. A \emph{direct system in }\index{direct system} $\C$ is a pair $\left (\{
    X_i\}_{i\in I}, \{\varphi_{i,j}\}_{i\leq j}\right )$ where $X_i$ is an object of $\C$ for every $i \in I$ and $\varphi_{i,j} \colon X_i \to X_j$ is a morphism in $\C$ for every $i \leq j$ such that
    \begin{enumerate}
        \item $\varphi_{i,i}$ is the identity,
        \item $\varphi_{i,k} = \varphi_{j,k}\circ \varphi_{i,j}$ for every $i \leq j \leq k$. 
    \end{enumerate}
    A pair $(X, \{\varphi_{i}\}_{i \in I})$, where $X$ is an object in $\C$ and $\varphi_i: X_i \to X$ is a morphism for every $i\in I$, is a \emph{cocone  of $(\{X_i\}_{i\in I}$, $\{\varphi_{i,j}\}_{i \leq j})$}\index{cocone } if $\varphi_i = \varphi_j \circ \varphi_{i,j}$. We say that the cocone  $(X,\{\varphi_i\}_{i\in I})$ is a \emph{direct limit}\index{direct!limit} if for any other cocone  $(Y, \{\psi_i\}_{i\in I})$ there is an unique morphism $\theta \colon X \to Y $ such that $\psi_i = \theta \circ \varphi_i$ for every $i \in I$. The object $X$ usually is denoted as $X = \varinjlim X_i$. This can be represented as the commutative diagram of Figure \ref{fig_limite_directo}.
    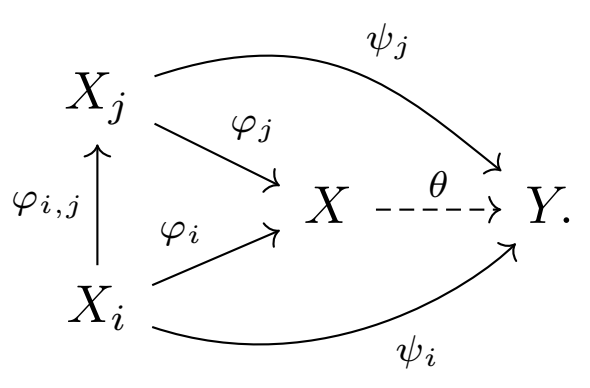
\begin{figure}[h]
        \centering
        \begin{tikzcd}[row sep=tiny]
         X_j\arrow[dr,"\varphi_j"] \arrow[drr,controls={+(1.5,0.5) and +(-1,0.8)}, "\psi_j"]&       &     \\
         & X \arrow[r,dashed,"\theta"] & Y.\\
         X_i \arrow[ur,"\varphi_i"] \arrow[uu,"\varphi_{i,j}"] \arrow[urr,controls={+(1.5,-0.5) and +(-0.5,-0.5) }, "\psi_i"']&       &   
        \end{tikzcd}
        \caption[Universal property of the direct limit]{Commutative diagram of the direct limit of a direct system. }
        \label{fig_limite_directo}
    \end{figure}
\end{definition}

\begin{definition}
Let $\mathcal{X}=\left(\{X_i\}_{i\in I},\{\varphi_{ij}\}_{i\leq j}\right)$ and $\mathcal{Y}=\left(\{Y_i\}_{i\in I},\{\psi_{ij}\}_{i\leq j}\right)$ be two direct systems indexed by the same directed set $I$.
An \emph{isomorphism of direct systems} from $\mathcal{X}$ to $\mathcal{Y}$ (denoted $\mathcal{X}\cong\mathcal{Y}$) is a family of isomorphisms $f_i\colon X_i\longrightarrow Y_i$  $i\in I$ such that, for all $i\leq j$, the following diagram commutes:
$$
\minCDarrowwidth 20pt
\begin{CD}
X_i @>{\varphi_{ij}}>> X_j\\
@V{f_i}VV @VV{f_j}V\\
Y_i @>{\psi_{ij}}>> Y_j.
\end{CD}
$$
\end{definition}
It is straightforward to realize that if two direct systems are isomorphic, then  $\varinjlim \mathcal{X}$ exists if and only if  $\varinjlim \mathcal{Y}$ exists also. In this case, there is an
isomorphism $f\colon\varinjlim\mathcal{X}\to\varinjlim\mathcal{Y}$ making commutative the squares:
\[
\minCDarrowwidth 15pt
\begin{CD}
X_i @>{\varphi_{i}}>> \varinjlim\mathcal{X}\\
@V{f_i}VV @VV{f}V\\
Y_i @>{\psi_{i}}>> \varinjlim\mathcal{Y}.
\end{CD}
\]

In general, not every direct system has a direct limit. If a category $\C$ has all arbitrary small coproducts and all coequalizers, any arbitrary colimit over a small diagram \(F: I \to \mathcal{C}\) can be explicitly constructed.

 For example, the category of groups, algebras, rings, modules over a ring, etc. We include some examples for self-containment, however, these are well known and can be checked in \cite{maclane}. In this vein,  recall that
if $R$ is a ring and $(\{M_i\}_{i\in I},\{\varphi_{i,j}\})$ a direct system of $R$-modules, then $\varinjlim M_i=\bigsqcup_{i\in I} M_i / \sim 
$ where $x_i \in M_i$ and $x_j \in M_j$ are related if there is some $k \geq i,j$ such that 
    $\varphi_{i,k} (x_i) = \varphi_{j,k}(x_j)$. The elements of $\varinjlim M_i$ are denoted $[x_i]$. The operations in the limit are 
    \begin{equation*}
    \begin{array}{rl}
         r \cdot [x_i] & := [r x_i] \ \forall r \in R, \cr
         [x_i]+[x_j] & := [\varphi_{i,k}(x_i)+\varphi_{j,k}(x_j)] \text{ with } i,j\leq k.
    \end{array}
    \end{equation*}
Limits in the category of $R$-algebras can be constructed as the limit in the category of $R$-modules, enriched with a product  
$[x_i][x_j] = [\varphi_{i,k}(x_i)\varphi_{j,k}(x_j)]$ with $i,j \leq k$.

The possibility of writing certain graphs as a suitable direct limit of finite graphs has been studied in \cite{AraMorenoPardo} and \cite{Goodearl}. Following \cite{Goodearl} we denote by $\Dg$ the category of directed graphs and graph homomorphisms, and by $\Dg_0$ that of directed graphs, whose homomorphisms $f\colon E\to F$ have injective $f^0\colon E^0\to F^0$. 
For a graph homomorphism $f\colon E\to F$ we denote by
$\tilde{f}\colon KE\to KF$ the unique $K$-algebra homomorphism whose restriction to $E^0\cup E^1$ is $f$.
By \cite[Lemma 2.5]{Goodearl}, (1) arbitrary direct limits exist in the categories $\Dg$ and $\Dg_0$; and (2) for any field $K$, the functor $K\_\colon\Dg_0\to\Alg$ such that $E\mapsto KE$ and $f\mapsto\tilde{f}$ preserves direct limits. Thus, if $E$ is a graph and $E=\varinjlim{S}$ where $S=(\{E_i\}_{i\in I},\{\varphi_{i,j}\}_{i\le j})$ is a direct system in $\Dg_0$, then $KE\cong\varinjlim{KS}$ where $KS=(\{KE_i\}_{i\in I},\{\widetilde{\varphi_{i,j}}\}_{i\le j})$. Also, we can consider the category $\CKD$ of arbitrary directed graphs and whose morphisms are arbitrary CK-morphisms \cite{Goodearl}  or complete morphisms in the terminology of \cite{AraMorenoPardo}. This is a subcategory of $\DG_0$, and the functor $L_K \colon \CKD \to \Alg$ (defined in \cite[\S 2.4]{Goodearl}) such that $E \mapsto L_K(E)$ preserves direct limits (\cite[Lemma 2.5]{Goodearl}). 

Given a graph $E$ and two finite subgraphs $F_i$ of $E$ ($i=1,2$), there is a finite subgraph of $E$ containing $F_i$ ($i=1,2$). Hence the family of finite subgraphs of $E$ is directed and so any graph $E$ is the union of its class of finite subgraphs. Thus, any graph is a direct limit in $\Dg_0$ of its finite subgraphs. Furthermore, if $E$ is a row finite graph, then it is the direct limit of its finite complete subgraphs.

\section{Limit of filtered algebras}\label{sec_limit_filtration}
In this section, we explore the concept of direct limit of a family of filtered algebras. In the category $\Falg$ this limit exists, however, in general, it is not in $\falg$. That is why we need to be cautious when we want to define a new notion of entropy. 

\begin{remark}
Consider the forgetful functor $U\colon\Falg\to\Alg$
acting on objects as $U(A,\F):=A$ (similarly for 
$V\colon\falg\to\Alg$). Both such functors are
faithful \cite{maclane}. If $S=\{(A_i,\F_i),\{\varphi_{i,j}\}_{i\le j}\}_{i\in I}$ is a direct system in $\Falg$, then define $U(S)$ as the direct system $U(S):=(\{A_i\}_{i\in I},\{\varphi_{i,j}\}_{i\le j})$ in $\Alg$. Then it is easy to check that if $\varinjlim S$ exists, we have
$\varinjlim S=(\varinjlim U(S),\{\varphi_i\}_{i\in I} )$ where the $\varphi_i\colon A_i\to \varinjlim U(S)$ are the canonical homomorphisms. On the other hand, if $\F_i=\{V_n^i\}_{n=1}^\infty$, then $(\{V_n^i\}_{i\in I},\{\varphi_{i,j}\vert_{V_n^i}\}_{i\le j})$ is a
direct system in $\Vect$.

Take a direct system of filtered algebras $\{(A_i,\F_i)\}_{i \in I}$ where $\F_i = \{V_n^i\}_{n=1}^\infty$ is a filtration and $\varphi_{i,j} \colon (A_i,\F_i) \to (A_j,\F_j)$ is a filtered algebra morphism for every $i \leq j$. Taking into account the above formula 
$\varinjlim S=(\varinjlim U(S),\{\varphi_i\}_{i\in I} )$
we have that the pairs $\Bigl (\{A_i\}_{i\in I},\{\varphi_{i,j}\}_{i \leq j}\Bigr )$ and $ \left (\{V_n^i\}_{i \in I}, \{\left .\varphi_{i,j}\right |_{V_n^i}\}_{i \leq j} \right )$ for any $n \geq 1$, are direct systems in $\Alg$ and $\Vect_K$ respectively.    
\end{remark}

Following the previous remark, it makes sense to consider the direct limit of the family of algebras once we forget that they are also filtered. This leads to the construction of the limit of filtered algebras.
\begin{proposition}\label{prop_limit_filteredalgebras}
   Let $S = \left ( \left \{(A_i,\F_i) \right \}_{i\in I}, \{\varphi_{i,j}\}_{i \leq j} \right )$ be a direct system of filtered algebras with $\F_i = \{V_n^ i\}_{n=1}^\infty$. If we consider the direct system $T=(\{A_i\}_{i\in I}, \{\varphi_{i,j}\}_{i\le j})$ in the category $\Alg$ whose limit is $\varinjlim T=(A,\{\varphi_i\}_{i\in I})$. Then defining $\F = \{V_n\}_{n=1}^\infty$ with $V_n = \bigcup_{i\in I} \varphi_i(V_n^i)$, one has that $\left ((A,\F), \{\varphi_i\}_{i\in I} \right )$ is the direct limit of the system $S$ (in $\Falg$).
\end{proposition}
\begin{proof}
    We need to check that $(A,\F) \in \Falg$, $\varphi_i$ is a morphism of filtered algebras 
 for every $i \in I$, $\left ((A,\F), \{\varphi_i\}_{i\in I} \right )$ is a cocone  and it verifies the universal property. 
Take $\hat{A} = \bigsqcup A_i /\sim$, since $\hat{A} \cong \varinjlim A_i$, by the universal property we can consider that $A = \hat{A}$ and $\varphi_i \colon A_i \to A$ with $\varphi_i(a) = [a] \in A$. To see that $\F$ is a filtration we follow Definition \ref{def_entropy_filtracion}:
    \begin{enumerate}
        \item By construction, we have that $V_n \subseteq A$ is a vector space for every $n\geq 1$ and they form a chain, that is, $V_1 \subseteq V_2 \subseteq \cdots \subseteq V_n \subseteq V_{n+1} \subseteq \cdots $. 
        \item If $[x] \in A$, then there is $i\in I$ such that $x \in A_i$. Since $\F_i$ is a filtration in $A_i$, there is some $n$ such that $x \in V_n^i$ and $\varphi_i(x) = [x] \in V_n$.
        \item Take $[x] \in V_n$ and $[y] \in V_m$, then there are some $i,j \in I$ such that $x \in V_n^i$ and $y \in V_m^j$. If we take $k \geq i,j$, then $\varphi_{i,k}(x) \in V_n^k$, $\varphi_{j,k}(y) \in V_m^k$ and $\varphi_{i,k}(x)\varphi_{j,k}(y) \in V_{n+m}^k$ which means that $[x][y] = [\varphi_{i,k}\varphi_{j,k}(y)]\in V_{n+m}$. 
    \end{enumerate}
We already know that $\varphi_i$ is a morphism of algebras and by definition it verifies that $\varphi_i(V_n^i) \subseteq V_n$. Consequently, $\varphi_i\colon (A_i,\F_i) \to (A,\F)$ is a filtered algebra morphism. In addition, since $(A,\{\varphi_{i}\}_{i\in I})$ is already a direct limit for $(\{A_i\}_{i\in I},\{\varphi_{i,j}\}_{i \leq j})$, it satisfies that $\varphi_i = \varphi_j \varphi_{i,j}$. 
Finally, we just need to check that $((A,\F),\{\varphi_i\}_{i\in I})$ satisfies the universal property. Take a cocone  $((B,\G), \{\psi_i\}_{i \in i})$, we need to find an unique $\theta \colon (A,\F) \to (B,\G)$ such that $\psi_i = \theta \circ \varphi_i$ for every $i \in I$. If we have $x_i \in A_i$, then $\theta$ has to verify that $\psi_i(x) = \theta(\varphi(x_i)) = \theta([x_i])$. In this sense we define $\theta\colon A \to B$ such that $\theta ([x]) := \psi_i(x_i)$ where $x_i \in A_i$ and $[x_i] = [x]$. This map is well-defined and is unique by construction. Indeed, take some $x_j \in A_j$ with $[x_j] = [x_i]$, then there is some $k \geq i,j$ such that $\varphi_{j,k}(x_j) = \varphi_{i,k}(x_i)$. Thus, if we compose with $\psi_k$ we have that 
\begin{equation*}
    \theta([x_j]) = \psi_j(x_j) = \psi_k(\varphi_{j,k}(x_j))= \psi_k( \varphi_{i,k}(x_i)) = \psi_i(x_i)  =\theta([x_i]).
\end{equation*}
Lastly, it is a filtered algebra morphism. Consider $\G = \{W_n\}_{n=1}^\infty$, if $[x] \in V_n$, then there is some $x_i \in V_n^i$ such that $[x_i] = [x]$ and $\theta([x]) = \psi_i(x_i) \in W_n$ because $\psi_i$ is a filtered algebra morphism. In addition, if $\lambda \in K$, $x_i \in A_i$, $y_j \in A_j$ and $k \geq i,j$, then 

\begin{gather*}
    \theta(\lambda [x_i]+[y_j]) = \theta([\varphi_{i,k}(\lambda x_i)+\varphi_{j,k}(x_j)]) = \psi_k(\varphi_{i,k}(\lambda x_i)+\varphi_{j,k}(x_j)) = \\
    \lambda \psi_k(\varphi_{i,k}(x_i)+\psi_k(\varphi_{j,k}(y_j))= \lambda \psi_i(x_i) +\psi_j(y_j) = \lambda\theta([x_i])+ \theta([y_j]).
\end{gather*}
and 
\begin{gather*}
    \theta([x_i][y_j]) = \theta([\varphi_{i,k}(x_i)\varphi_{j,k}(x_j)]) = \psi_k(\varphi_{i,k}(x_i)\varphi_{j,k}(x_j)) = \\
    \psi_k(\varphi_{i,k}(x_i))\psi_k(\varphi_{j,k}(y_j))= \psi_i(x_i) \psi_j(y_j) = \theta([x_i]) \theta([y_j]).
\end{gather*}
\end{proof}
\begin{remark}
    The reader can check that in the context of the previous Proposition, $V_n = \varinjlim V_n^i$ since $(V_n,\{\left .\varphi_i\right |_{V_n^i} \}_{i\in I})$ satisfies the universal property. To check so we only need to find an unique map $\theta \colon V_n \to \bigsqcup V_n^i / \sim $ such that $\phi_i = \theta \circ \varphi_i$ with $\phi_i \colon V_n^i \to V_n^i / \sim $ the map that sends $x$ to $\overline{x} \in V_n^i / \sim $. We change the notation for the new equivalence class in order to distinguish the one in $V_n$ to the one in $V_n^i / \sim $. The needed $\theta$ is just the map that satisfies $\theta([x_i]) = \overline{x_i}$. It is well-defined since $[x_i] = [x_j]$ if and only if there is some $k \geq i,j$ such that $\varphi_{i,k}(x_i) = \varphi_{j,k}(x_j)$ which is equivalent to $\overline{x_i} = \overline{x_j}$. Finally, highlight that since $V_n = \varinjlim V_n^i$ we will not make distinction between them or their respective equivalence classes and use the $[\cdot]$ notation. 
\end{remark}

    Even though the direct limit exists in the category of filtered algebras, this is not true in the categories of finitely filtered algebras. More specifically, the direct limit of finitely filtered algebras, as in Proposition \ref{prop_limit_filteredalgebras},  is not (in general) a finitely filtered algebra. 

    The path algebra $KCh_\infty$ with $Ch_\infty$ the infinity chain is a good example of this. Consider the direct system $(\{(KCh_i,\F_i)\}_{i=2}^\infty,\{\varphi_{i,j}\}_{i,j=2}^\infty)$ with $\varphi_{i,j}\colon KCh_i\hookrightarrow KCh_j$ the natural inclusion for $j \geq i$ and $\F_i$ the standard filtration of $KCh_i$. In this case, the direct limit is $\left ((KCh_\infty,\F),\{\varphi_i\}_{i=2}^\infty \right )$ with $\varphi_i \colon KCh_i \to KCh_\infty$ the inclusion and $\F$ the standard filtration (a way to prove this is using \cite[Lemma 2.5]{Goodearl}). However, as we will see, the standard filtration is not finite for the path algebra of the infinity chain whereas the standard filtrations for the graphs in the direct system are. 

    Instead of a particular example, we can see more generally that the limit of finitely filtered algebras is not (in general) a finitely filtered algebra. We can see this by computing the dimension of the quotient of the limit. We will need for that the next results.
    Some of them can be considered general knowledge and a person who is familiar with category theory and direct limits might think of them as a triviality. If that is the case, the reader can skip the proofs of the next results until Section \ref{sec_entropy_directsystems}. However, since we want to be as self-contained as we can and we do not want to deviate the content from our problem, we include the explicit proofs of the desired lemmas. 
    
\begin{lemma}\label{lemma_entropy_limite_injectivo}
    If $S = (\{A_i\}_{i\in I}, \{\varphi_{i,j}\}_{i,j\in I} )$ is a direct system where $\varphi_{i,j}\colon A_i \to A_j$ is a monomorphism for every $i \leq j$, then the morphisms in the direct limit $(A,\{\varphi_i\}_{i\in I})$ of $S$ are injective for every $i\in I$. Furthermore, there is an isomorphic system $S' = (\{A'_i\}_{i\in I}, \{\varphi'_{i,j}\}_{i,j\in I} )$ with limit $(A,\{\varphi'_i\}_{i\in I})$ where each $A_i'$ is a subalgebra 
    of $A$, and both, the $\varphi'_{i,j}$ and the $\varphi'_i$ are inclusions.
\end{lemma}
\begin{proof}
    First of all, consider $\hat{A} = \bigsqcup A_i /\sim$, since $\hat{A}$ also is the direct limit of $S$, then $A \cong \hat{A}$. By the universal property, we just need to prove that the morphisms $A_i \to \hat{A}$ are injective. Remember that 
    \begin{equation*}
        \begin{array}{cccc}
             \psi_i \colon & A_i& \to & \hat{A}  \\
             & x&\mapsto & [x] 
        \end{array}
    \end{equation*}
    and consider $x,y \in A_i$ with $[x] = [y]$. This means that there is some $k \geq i$ such that $\varphi_{i,k}(x) = \varphi_{i,k}(y)$ and since all the morphisms $\varphi_{i,j}$ are injective, $x = y$ as we expected. For the second assertion, define for each $i$ the algebra $A_i':=\varphi_i(A_i)$ and for $i\le j$ notice that $A_i'\subseteq A_j'$. Then define $\varphi'_{i,j}$ as the inclusión $A_i'\hookrightarrow A_j'$. The direct system so constructed $S'$ is isomorphic to $S$.
\end{proof}
Consider $S = (\{(A_i,\F_i)\}_{i\in I},\{\varphi_{i,j}\}_{i,j\in I})$ a direct system of filtered algebras with $\F_i = \{V_n^i\}_{n=1}^\infty$ for $i \in I$. The map $\varphi_{i,j}$ is a filtered algebra morphism, then $\varphi_{i,j}(V_{n-1}^i)\subseteq V_{n-1}^j$. In a similar way to what we did for \cite[Lemma 4.1]{entropy1} and \cite[Lemma 4.2]{entropy1} , we can induce a vector space morphism $\bar{\varphi}_{i,j} \colon 
V_n^i/V_{n-1}^i \to V_n^j/V_{n-1}^j$ such that $\bar{\varphi}_{i,j}(x+V_{n-1}^i) = \varphi_{i,j}(x) + V_{n-1}^j$.

\begin{lemma}\label{lemma_entropy_limitedelcociente_cocientedellimite}
    Let $S$ be a direct system of filtered algebras as before. The family of vector spaces $\{V_n^i/V_{n-1}^i\}_{i\in I}$ together with the morphisms $\{\bar{\varphi}_{i,j}\}_{i \leq j}$ is a direct system and for every $n \geq 2$
   \begin{equation}\label{eq_entropy_isomorfo}
   \frac{\bigsqcup_{i \in I}V_n^i/\sim}{\bigsqcup_{i \in I}V_{n-1}^i  /\sim } \cong \frac{\bigsqcup_{i\in I} V_n^i/V_{n-1}^i}{\sim }
   \end{equation}
   as $K$-vector spaces. 
\end{lemma}
\begin{proof}
First we need to check that $(\{V_n^i/V_{n-1}^i\}_{i\in I},\{\bar{\varphi}_{i,j}\}_{i \leq j} )$ is a direct system. If $x_i + V_{n-1}^i \in V_n^i /V_{n-1}^i$ and $i \leq j \leq k$, then
\begin{equation*}
    \bar{\varphi}_{i,i}(x_i +V_{n-1}^i) = \bar{\Id}(x_i+V_{n-1}^i) = x_i + V_{n-1}^i,
\end{equation*}
and
\begin{gather*}
    \bar{\varphi}_{i,k}(x_i + V_{n-1}^i) = \varphi_{i,k}(x_i) +V_{n-1}^k = \varphi_{j,k}(\varphi_{i,j}(x_i)) + V_{n-1}^k = \\
    \bar{\varphi}_{j,k}(\varphi_{i,j}(x_i) + V_{n-1}^j) = \bar{\varphi}_{j,k}\circ \bar{\varphi}_{i,j}(x_i + V_{n-1}^i).
\end{gather*}
To prove \eqref{eq_entropy_isomorfo} consider $V_n =\bigsqcup_{i \in I}V_n^i/\sim$ and $V_{n-1} = \bigsqcup_{i \in I}V_{n-1}^i/\sim$. We define the map 
    \begin{equation*}
    \begin{array}{cccc}
         \Omega \colon  & \frac{\bigsqcup_{i\in I} V_n^i/V_{n-1}^i}{\sim } & \to &  V_n/V_{n-1} \\
         & [x_i + V_{n-1}^i] & \mapsto &  [x_i]+ V_{n-1}
    \end{array}
    \end{equation*}
    which is an isomorphism of $K$-vector spaces. In order to prove so we need to verify that it is well-defined, linear and bijective:
    \begin{enumerate}
        \item Consider two elements $x_i + V_{n-1}^i $ and $x_j + V_{n-1}^j$ related. This means that there is certain $k \in I$ with $k \geq i,j$ such that $\bar{\varphi}_{i,k}(x_i + V_{n-1}^i) = \bar{\varphi}_{j,k}(x_j + V_{n-1}^j)$. In particular, $\varphi_{i,k}(x_i) + V_{n-1}^k = \varphi_{j,k}(x_j)+V_{n-1}^k$ which implies that $\varphi_{i,k}(x_i) - \varphi_{j,k}(x_j) \in V_{n-1}^k$. Then, $[\varphi_{i,k}(x_i) - \varphi_{j,k}(x_j)] = [x_i]- [x_j] \in V_{n-1}$ and
        \begin{equation*}
            \Omega\left ([x_i+V_{n-1}^i] \right ) = [x_i]+ V_{n-1} = [x_j]+V_{n-1} = \Omega \left ([x_j + V_{n-1}^j] \right ).
        \end{equation*}
        \item Take $\lambda \in K$, then 
        \begin{gather*}
            \Omega(\lambda \cdot [x_i + V_{n-1}^i] ) = \Omega ([\lambda x_i + V_{n-1}^i]) = [\lambda x_i] + V_{n-1} = \\
            \lambda [x_i] + V_{n-1} = \lambda \cdot \left ([x_i] + V_{n-1} \right ) = \lambda \Omega([x_i + V_{n-1}^i]).
        \end{gather*}
        Consider $[x_i+V_{n-1}^i],[x_j+V_{n-1}^j] \in \frac{\bigsqcup_{i\in I} V_n^i/V_{n-1}^i}{\sim }$ then we have
        \begin{gather*}
            \Omega\left ([x_i+V_{n-1}^i]+[x_j+V_{n-1}^j]\right )= \Omega\left ([\bar{\varphi}_{i,k}(x_i+V_{n-1}^i)+ \bar{\varphi}_{j,k}(x_j+V_{n-1}^j)] \right ) = \\
            = \Omega\left ([\varphi_{i,k}(x_i)+\varphi_{j,k}(x_j)+V_{n-1}^k]\right ) = [\varphi_{i,k}(x_i)+\varphi_{j,k}(x_j)]+V_{n-1} = \\
            = [x_i]+ [x_j] + V_{n-1} = \left ([x_i]+ V_{n-1}\right )+\left ([x_j]+ V_{n-1} \right) = \\
            = \Omega\left ([x_i+V_{n-1}^i]\right )+\Omega\left ([x_j+V_{n-1}^j] \right).
        \end{gather*}
        \item To prove the injectivity, take $[x_i+ V_{n-1}^i ]$ with $\Omega([x_i+ V_{n-1}^i ]) = 0$. That is $[x_i] + V_{n-1} = 0$, which means that $[x_i] \in V_{n-1}$. In particular, there is a $x_j \in V_{n-1}^j$ such that $[x_j] = [x_i] \in V_{n-1}$, that is, there is $k \geq i,j$ with $\varphi_{i,k}(x_i) = \varphi_{j,k}(x_j)$. Since, $\varphi_{j,k}$ is a filtered algebra morphism and $x_j \in V_{n-1}^j$ then $\varphi_{i,k}(x_i) = \varphi_{j,k}(x_j) \in V_{n-1}^k$. Then, $\varphi_{i,k}(x_i) + V_{n-1}^k = 0$ and $\bar{\varphi}_{i,k}(x_i + V_{n-1}^i) = 0$. Consequently, by the relation $\sim$ in $\frac{\bigsqcup_{i\in I} V_n^i/V_{n-1}^i}{\sim }$, this means that $[x_i+V_{n-1}^i] = 0$. The surjectivity of the map is immediate by construction. If we take $[x_i] + V_{n-1} \in V_n/V_{n-1}$, then $\Omega([x_i+V_{n-1}^i]) = [x_i] + V_{n-1}$.
    \end{enumerate}
\end{proof}

\begin{remark}\label{rem_}
There is a smarter way to prove this using exact functors. If $\Vect_K$ denotes the category of $K$-vector spaces and $I$ is directed set, then $\Vect_K^I$ is the category which objects are direct systems of vector spaces. Its objects are direct systems $(\{V_i\}_{i\in I},\{\varphi_{i,j}\}_{i\leq j})$ and if $S,S' \in \Vect_K^I$ then a morphism $\theta \colon S \to S'$ is a family of maps $\{\theta_i\}_{i\in I}$ with $\theta_i\colon V_i \to V_i'$ such that $\varphi_{i,j}' \circ \theta_i = \theta_j \circ \varphi_{i,j}$ for every $i\leq j$. Under this context, we can define the functor $\F\colon \Vect_K^I \to \Vect_K$ such that if $S =(\{V_i\}_{i\in I},\{\varphi_{i,j}\}_{i\leq j})$, $S' = (\{V'_i\}_{i\in I},\{\varphi'_{i,j}\}_{i\leq j})$ and $\theta \colon S \to S'$, then $\F(S) := \varinjlim V_i$ and $\F(\theta) \colon \varinjlim V_i \to \varinjlim V'_i$ is the only morphism given by the universal property such that $\varphi_i'\circ \theta_i  = \F(\theta) \circ \varphi_i$. This functor is exact (\cite[Theorem 2.6.15]{weibel1994introduction}), that is, if 
\begin{equation*}
    0 \to  S_1 \hookrightarrow S_2 \twoheadrightarrow S_3 \to 0
\end{equation*}
is an exact sequence, then 
\begin{equation*}
    0 \to \F(S_1) \hookrightarrow \F(S_2) \twoheadrightarrow \F(S_3) \to 0
\end{equation*}
is also an exact sequence. In particular, if $(\{(A_i,\F_i)\}_{i\in I},\{\varphi_{i,j}\}_{i,j\in I})$ is a direct system of filtered algebras with $\F_i = \{V_n^i\}_{n=1}^\infty$ then there is a exact sequence in $\Vect_K^I$ induced by 
\begin{equation*}
    0 \to V_{n-1}^i \hookrightarrow V_{n}^i \twoheadrightarrow V_n^i/V_{n-1}^i \to 0 .
\end{equation*}
Consequently,
\begin{equation*}
0 \to \varinjlim V_{n-1}^i \hookrightarrow \varinjlim V_{n}^i \twoheadrightarrow \varinjlim V_n^i/V_{n-1}^i \to 0 
\end{equation*}
is an exact sequence because $\F$ is an exact functor. Finally this implies that 
\begin{equation*}
    \varinjlim V_n^i/\varinjlim V_{n-1}^i \cong \varinjlim V_n^i/V_{n-1}^i.
\end{equation*}
\end{remark}

\begin{lemma}\label{lemma_entropy_limit_dimension}
   Consider $(\{V_i\}_{i\in I},\{\varphi_{i,j}\}_{i\leq j})$ a direct system of vector spaces with $\varphi_{i,j}\colon V_i \to V_j$ injective for every $i\leq j$. If $V = \varinjlim V_i$, then 
    \begin{equation}
        \dim V = \sup_{i\in I}\{\dim V_i\}
    \end{equation}
\end{lemma}
\begin{proof}
    Firstly, we consider the case that $\dim V = n < +\infty$. By Lemma \ref{lemma_entropy_limite_injectivo} we have that $\dim V_i \leq \dim V$ for every $i \in I$. Take a basis $\mathcal{B} = \{[x_{i_1}],[x_{i_2}],\ldots,[x_{i_n}]\}$. Since $I$ is directed, there is some $k\in I$ such that $k \geq i_j$ for every $j=1,\ldots, n$. In addition, $[x_{i_j}] = [\varphi_{i_j,k}(x_{i_j})]$ for every $j=1,\ldots,n$ and without loss of generality we can consider a basis $\mathcal{B} = \{[x_1],[x_2],\ldots, [x_n]\}$ of $V$, with $x_i \in V_k$ for every $i = 1,\ldots, n$. Since $B$ is a set of linearly independent vectors in $V$, it is easy to check that $\{x_1,\ldots,x_n\}$ are linearly independent in $V_k$. Consequently $n \leq \dim V_k \leq n $  and $\dim V = \sup_{i\in I}\{\dim V_i\}$. 

    Secondly, if $\dim V = +\infty$, then there is a sequence of finite-dimensional vector spaces $\{W_n\}_{n=1}^\infty$ such that $W_n \subseteq V$ and $\lim_{n \to \infty} \dim W_n = +\infty$. Using the previous argument, we can construct a sequence of vector spaces $\{V_{i_n}\}_{n=1}^\infty$ such that $\dim W_n \leq \dim V_{i_n}$. Consequently,
    \begin{equation*}
        \dim V = +\infty = \sup_{i\in I}\{\dim V_i\}.
    \end{equation*}
\end{proof}

\begin{proposition}\label{prop_limitefinito_nofinito}
   The category $\falg$ is not closed for direct limits. 
\end{proposition}
\begin{proof}
    Consider the path algebra $KCh_\infty$. Following \cite[Lemma 3.1]{AraMorenoPardo} and \cite[Lemma 3.5]{Goodearl}, we have that $KCh_{\infty}=\varinjlim S$, where  $S= (\{KCh_i\}_{i \geq 2}, \{\varphi_{i,j}\}_{i\leq j})$ and the transition homomorphisms are $\varphi_{i,j}\colon KCh_i \to KCh_j$, the algebra morphisms induced by the inclusions $Ch_i \hookrightarrow Ch_j$  of finite graphs. For each $i \geq 2$, we can consider the standard filtration of path algebras of finite graph $\F_i$. The inclusion of graphs is compatible with the standard filtration of finite graphs, that is, the inclusion of algebras is a strict monomorphism of filtered algebras. Under these conditions, we define the strict direct system of finitely filtered algebras $\hat{S} = (\{(KCh_i,\F_i)\}_{i\geq 2},\{\varphi_{i,j}\}_{i\leq i})$. By Proposition \ref{prop_limit_filteredalgebras} we have that $(KCh_\infty,\F)$ with $\F = \{V_n\}_{n\geq 1}$ and $V_n = \varinjlim V_n^i$ is the direct limit of $\hat{S}$ in $\Falg$. However, $(KCh_\infty,\F)$ is not in $\falg$ that is, $\F$ is not a finite filtration. On the line of Lemmas \ref{lemma_entropy_limitedelcociente_cocientedellimite} and \ref{lemma_entropy_limit_dimension}, we have that 
\begin{equation}\label{sinast}
    \dim V_n/V_{n-1} = \dim \varinjlim V_n^i/V_{n-1}^i = \sup_{i\in I} \left \{\dim V_n^i/V_{n-1}^i \right \}.
\end{equation}
Since we are working with the standard filtration and by \cite{entropy1} we have that
\begin{equation*}
    \dim V_n/V_{n-1} = \sup_{i\geq 2}\{\left |\lambda \in \PP(Ch_i) \colon |\lambda| = n \right |\} = +\infty
\end{equation*}
\end{proof}
\section{Entropy of strict direct systems}\label{sec_entropy_directsystems}
It is important not to lose the spotlight in our problem. We are trying to obtain a definition of the entropy that generalise the one that we have for filtered algebras and that is compatible with the previous one. To do so we take inspiration from the infinity chain. In this case, for a field $K$, we have that $KCh_i$ is included in $KCh_\infty$ by the usual inclusion $\mathfrak{i}\colon KCh_i \hookrightarrow KCh_\infty$. In addition, if $\F=\{V_n\}_{n=1}^\infty$ is the standard filtration of $KCh_\infty $ and $\F_i = \{V_n^i\}_{n=1}^\infty$ is the standard filtration of $KCh_i$, then $\mathfrak{i}^{-1}(V_n) = V_n^i$. Following these ideas, we define:
\begin{definition}\label{def_entropy_strict}
    Let $(A,\F)$ and $(B,\G)$ be two filtered algebras with $\F = \{V_n\}_{n=1}^\infty$ and $\G = \{W_n\}_{n=1}^\infty$. We say that a monomorphism of filtered algebras $f \colon (A,\F) \to (B,\G)$ is \emph{strict}\index{monomorphism!strict} if and only if $V_n = f^{-1} (W_n)$.
\end{definition}
\begin{definition}\label{def_entropy_strict_system}
    A direct system of filtered algebras $S = \left (\{(A_i,\F_i)\}_{i \in I}, \{\varphi_{i,j}\}_{i \leq j}\right)$ is \emph{strict}\index{direct system!strict} if $\varphi_{i,j}$ is a strict monomorphism for every $i \leq j$.  
\end{definition}
\begin{lemma}\label{lemma_strict_morphism}
    If $S = \left (\{(A_i,\F_i)\}_{i \in I}, \{\varphi_{i,j}\}_{i \leq j}\right)$ is a strict direct system of filtered algebras with $\varinjlim S = ((A,\F), \{\varphi_i\}_{i \in I})$, then every $\varphi_i$ is a strict monomorphism of filtered algebras.
    Moreover, there is a family of subalgebras $\{A_i'\}_{i\in I}$ where $A_i\cong A_i'$ such that 
    $A=\bigcup_{i\in I} A_i'$ and $A_i'\subset A_j'$ when $i\le j$. Each $A_i'$ is a filtered algebra relative to $\F_i':=\{V_n\cap A_i'\}_{n\ge 1} = \{\varphi_i(V_n^i)\}_{n\geq1}$ (where $\F=\{V_n\}_{n\ge 1}$).
\end{lemma}
\begin{proof}
    By Lemma \ref{lemma_entropy_limite_injectivo} we have that  $\varphi_i \colon A_i \to A$ is a filtered algebra monomorphism. Lastly, we have to check that it is strict, that is $\varphi_i^{-1}(V_n) \subseteq V_n^i$. Take $a \in \varphi_i^{-1}(V_n) \subseteq A_i$, then $\varphi_i(a) \in V_n$, which means that there is and $j \in I$  and $a' \in V_n^j$ such that $\varphi_i(a) = \varphi_j(a')$. Since $I$ is a directed set, there is $k \geq i,j$, then $\varphi_i(a) = \varphi_k(\varphi_{i,k}(a)) = \varphi_k(\varphi_{j,k}(a')) $ and by injectivity of $\varphi_k$, we have $\varphi_{i,k}(a) = \varphi_{j,k}(a') \in \varphi_{j,k}(V_n^j) \subseteq V_n^k$. Since $\varphi_{i,j}$ are strict, then $a \in V_n^i$. The second part induced directly from the definition of strict monomorphism of filtered algebras and the properties of the strict direct system. That is $\varphi_j \varphi_{i,j} = \mathfrak{i}_{i,j}\varphi_i$ where $\mathfrak{i}_{i,j}$ is the inclusion of $\varphi_i(A_i)= A_i'$ in $\varphi_j(A_j) = A_j'$, in other words, the diagram
\begin{equation*}
    \begin{CD}
V_n^i @>{\varphi_{ij}}>> V_n^j\\
@V{\varphi_i}VV @VV{\varphi_j}V\\
\varphi_i(V_n^i) @>{\mathfrak{i}_{ij}}>> \varphi_j(V_n^j),
\end{CD}
\end{equation*}
is commutative, whence the stated result holds.
\end{proof}

Once we know the concepts of strict monomorphism and strict direct system, we can define the entropy of an algebra as the entropy of a strict direct system. 

\begin{definition}\label{def_entropy_limit}
    Let $A$ be a $K$-algebra such that there is a strict direct system $S = \left (\{(A_i,\F_i)\}_{i \in I}, \{\varphi_{i,j}\}_{i \leq j}\right)$ of finitely filtered algebras with $A = \varinjlim A_i$ as $K$-algebras. Then, we define the entropy of $A$ for the system $S$, denoted by $\h(A,S)$, as
    \begin{equation*}
        \h(A,S) := \sup_{i\in I}\{\h(A_i,\F_i)\}
    \end{equation*}
\end{definition}
Let us compute the entropy of the path algebra of the infinity chain.

\begin{example}\label{ejemplo_chorizo}
    Consider $S =(\{(KCh_i,\F_i)\}_{i \geq 1},\{\varphi_{i,j}\}_{i \leq j})$ a strict direct system of filtered algebras with $\F_i$ the standard filtration and $\varphi_{i,j}\colon KCh_i \hookrightarrow KCh_j$ the inclusion. The adjacency matrix of $Ch_i$ is $A_i = \sum_{j=1}^{i} (E_{j,j+1}+E_{j+1,j})$ where $E_{jk}$ is the elementary $(j,k)$-matrix of size $i+1$ (which unique nonzero entry is $1$ in the $(j,k)$ place). This matrix is symmetric. By \cite[Lemma 4.1]{entropy1} the sequence of entropies 
$h_i:=\h(KCh_i,\F_i)$ is non-decreasing $h_i\le h_{i+1}$ and bounded since each $Ch_i$ is a subgraph of the extended graph of a cycle with $i+1$ vertices, denoted by $\hat{C}_{i+1}$ (see Figure \ref{figg_ciclo_n}).
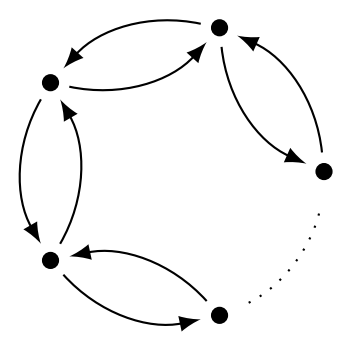
\begin{figure}[h]
    \centering
        \begin{tikzpicture}[shorten <=2pt,shorten >=2pt,>=latex, node distance={15mm}, main/.style = {draw, fill, circle, inner sep = 1pt}]

\def \n {5}
\def \radius {1cm}
\def \margin {8} 

\foreach \s in {1,...,4}
{
  \node[main] (\s) at ({360/\n * (\s - 1)}:\radius) {} ;
}
\node[main] (\n) at ({360/\n * (\n - 1)}:\radius) {};
  \draw[dotted] ({360/\n * (\n - 1)+\margin}:\radius) 
    arc ({360/\n * (\n - 1)+\margin}:{360/\n * (\n)-\margin}:\radius);

\draw[->] (1) to [bend right = 30] (2);
\draw[->] (2) to [bend right = 30] (1);
\draw[->] (2) to [bend right = 30] (3);
\draw[->] (3) to [bend right = 30] (2);
\draw[->] (3) to [bend right = 30] (4);
\draw[->] (4) to [bend right = 30] (3);
\draw[->] (4) to [bend right = 30] (5);
\draw[->] (5) to [bend right = 30] (4);

\end{tikzpicture}
    \caption{Extended directed cycle with $i$ vertices.}
    \label{figg_ciclo_n}
\end{figure}

The entropy of the  $K\hat{C}_{i+1}$ is $\log 2$ (\cite[Example 5.9]{entropy1}) and does not depend on the number of vertices, then we have an upper bound for $\{h_i\}_{i=2}^\infty$ which is $\log 2$. So $\{h_i\}_{i=1}^\infty$ is convergent and $\sup_{i \geq 2} \{h_i\}\le \log 2$. 

Our next goal is to prove the equality 
$\sup_{i \geq 2} \{h_i\}=\log 2$. From \cite[Theorem 5.4]{entropy1} we have that $h_i = \log \rho (A_i)$. Computing the spectral radius of the different matrices is a complex computational problem. However, we can obtain an smart lower bound that will prove our assertion. For every symmetric  matrix $A$ of size $m$ with $\spec(A)$ the set of its eigenvalues, we have that 
\begin{equation*}
    (A)_{1,1} \leq \tr(A) = \sum_{\lambda \in \spec(A)}\lambda \leq \sum_{\lambda \in \spec(A)}|\lambda| \leq m \rho(A).
\end{equation*}
Thus, the following inequality holds
\begin{equation*}
    (A^n)_{1,1} \leq m \rho(A^n) \leq m\rho(A)^n.
\end{equation*}
In particular $\sqrt[n]{(A^n)_{1,1}} \leq \sqrt[n]{i+1}\rho(A)$. If we apply this to the adjacency matrix $A_i$ of $Ch_i$ we have that
\begin{equation*}
    \sqrt[n]{(A_i^n)_{1,1}} \leq \sqrt[n]{i+1} \rho(A_i)
\end{equation*}
 for every $n \geq 1$. On the other hand, for every $n$, each entry $(A_i^{2n})_{1,1}$ counts the number of paths of length $2n$ from the vertex $v_1$ to $v_1$. It can be seen that this number is also the number of Dyck's words of length $2n$. By Corollary 6.2.3 (ii) or (iii) of \cite{stanley}, the number of Dyck's  words of length $2n$ is given by the Catalan number $c_n= \binom{2n}{n}\frac{1}{n+1}$. Then $(A_i^{2n})_{1,1} = c_n$ where $\{c_n\}_{n\geq 1}$ denotes the Catalan numbers. 
Thus, $\sqrt[2n]{c_n} \leq \sqrt[2n]{i+1} \rho(A_i)$ for every $n$ and every $i \geq 2$. Since the logarithm is an increasing function, we have that 
\begin{equation*}
    \frac{1}{2n}\log c_n \leq \frac{\log(i+1)}{2n}+\log \rho(A_i) = \frac{\log(i+1)}{2n}+h_i
\end{equation*}
for every $n$ and every $i\geq 2$. Using Stirling formula we have that $\lim_{n\to \infty} \frac{1}{2n}\log c_n = \log 2$.  Consequently, if $n$ tends to infinity, then $\log 2 \leq h_i$ for every $i \geq 2$. Finally,
\begin{equation*}
    \log 2 \leq \sup_{i\geq 2} h_i \leq \log 2
\end{equation*}
which means that the equality holds. In particular, $\h(KCh_\infty,S) = \log 2$.
\end{example}

Similarly as it happens for the entropy of filtered algebras, the entropy of an algebra which is direct limit of a strict system $S$ depends on the chosen system.
\begin{proposition}\label{prop_entropy_limit_dependencia}
    Consider $S = \left (\{(A_i,\F_i)\}_{i \in I}, \{\varphi_{i,j}\}_{i \leq j}\right)$ a strict direct system of finitely filtered algebras with $A = \varinjlim A_i$ and $\F_i= \{V_n^i\}_{i\in I}$. If we take $\G_i := \{W_n^i\}_{n=1}^\infty$ for a fixed $k \geq 1$ such that $W_n^i:= V_{kn}^i$, then we have the following:
    \begin{enumerate}
        \item $S' =\left (\{(A_i,\G_i)\}_{i \in I}, \{\varphi_{i,j}\}_{i \leq j}\right)$ is a strict direct system of finitely filtered algebras,
        \item $\h(A,S') = k \h(A,S)$.
    \end{enumerate}
\end{proposition}
\begin{proof}
First of all, we have to prove that $S'$ is a strict direct system of filtered algebras. By \cite[Proposition]{entropy1}, the family $\G_i$ is a filtration. In addition, by construction, the system is strict. The morphisms $\varphi_{i,j}$ remain being injective and $\varphi_{i,j}^{-1}(W_n^j) = \varphi_{i,j}^{-1}(V_{kn}^j) = V_{kn}^i = W_{n}^i$.
For the second part, observe that, by Proposition 
\begin{equation*}
    \h(A_i,\G_i) = k\h(A_i,\F_i).
\end{equation*}
Consequently,
\begin{gather*}
      \h(A,S') = \sup_{i\in I}\{\h(A_i,\G_i)\} =\\
      = \sup_{i \in I}\{k\h(A_i,\F_i)\} = k \sup_{i \in I}\{h(A_i\F_i)\} = k \h(A,S).
\end{gather*}
\end{proof}

\subsection{The Compatibility Question}\label{sub_compatibility}

Under this setting we ask ourselves if the entropy of filtered algebras and the entropy of direct limits are compatible. That is, if $(A,\F)$ is a finitely filtered algebra and $A = \varinjlim A_i$ with $(\{A_i\}_{i\in I},\{\varphi_{i,j}\}_{i,j\in I})$ a direct system where $\varphi_{i,j}\colon A_i \to A_j$ is a monomorphism, is there a strict direct system $S$ such that $\h(A,\F) = \h(A,S)$? The answer is positive, as we prove in Theorem \ref{th_entropy_limit_filtered}.

\begin{theorem}
\label{th_entropy_limit_filtered}
    Let $(A,\F)$ be a finitely filtered algebra with $\varinjlim A_i = (A, \{\varphi_i\}_{i\in I})$ for a direct system $(\{A_i\}_{i\in I},\{\varphi_{i,j}\}_{i,j\in I})$ where $\varphi_{i,j}\colon A_i \to A_j$ is a monomorphism for every $i\leq j$. The pair $S = \left (\{(A_i,\F_i)\}_{i\in I}, \{\varphi_{i,j}\}_{i\leq j}\right ) $ with $\F_i = \{\varphi_i^{-1}(V_n)\}_{n=1}^\infty$ is a strict direct system of finitely filtered algebras such that 
    \begin{equation*}
        \h(A,\F) = \h(A,S). 
    \end{equation*}
\end{theorem}
\begin{proof}
If $(A,\{\varphi_i\}_{i\in I})$ is the direct limit, then by Lemma \ref{lemma_entropy_limite_injectivo}, the map $\varphi_i$ is a monomorphism for every $i\in I$. Following Proposition \cite[Lemma 4.2]{entropy1} each $\F_i$ is a filtration and 
    \begin{equation*}
        \h(A_i,\F_i) \leq \h(A,\F)
    \end{equation*}
    for every $i \in I$. If we take the supremum, then
    \begin{equation*}
        \h(A,S) = \sup_{i\in I}\{\h(A_i,\F_i)\} \leq \h(A,\F).
    \end{equation*}
In order to prove that $\h(A,\F) \leq \h(A,S)$ we need to check that $\varinjlim V_n^i = V_n$.

First observe that $S$ is a strict direct system of finitely filtered algebras. Indeed, we only have to check that $\varphi_{i,j}^{-1}(V_n^j)=  V_n^i$. Since $\varphi_i = \varphi_j \circ \varphi_{i,j}$, then
    \begin{gather*}
        \varphi_{i,j}^{-1}(V_n^j) = \varphi_{i,j}^{-1}(\varphi_j^{-1}(V_n)) =  (\varphi_{j}\circ \varphi_{i,j})^{-1}(V_n) = \\
        = \varphi_i^{-1} (V_n) = V_n^j,
    \end{gather*}
    Consequently, by Lemma \ref{lemma_entropy_limitedelcociente_cocientedellimite} the pair $S_n = (\{V_n ^i\}_{i\in I}, \{\varphi_{i,j}\}_{i,j})$ is a direct system too. Finally, to verify that $\varinjlim V_n^i = V_n$ we are going to see that it verifies its universal property. Consider a cocone  $(W,\{\psi_i\}_{i \in I})$ for the direct system $S_n$. For a element $v \in V_n$, there exists a $i \in I$ such that $\varphi_i^{-1}(v) \in V_n^i$ (which is unique since $\varphi_i$ is a monomorphism). Then, we define the map
    \begin{equation*}
        \begin{array}{cccc}
             \theta \colon & V_n & \to & W   \\
             & v & \mapsto & \psi_i(\varphi_i^{-1}(v)).
        \end{array}
    \end{equation*}
    By construction it is direct that $\psi_i = \theta \circ \varphi_i$. However, we need to prove that the map $\theta$ is well-defined, that is, it does not depend on the chosen $i \in I$. Suppose $j \in I$ with $\varphi_{j}^{-1}(v) \in V_n^j$, there is some $k \geq i,j$ with  $\psi_i = \psi_k \circ \varphi_{i,k}$, $\psi_j = \psi_k\circ \varphi_{j,k}$, $\varphi_i = \varphi_k \circ \varphi_{i,k}$ and $\varphi_j = \varphi_k \circ \varphi_{j,k}$, then 
   \begin{gather*}
       \psi_i(\varphi_i^{-1}(v)) = \psi_k(\varphi_{i,k}((\varphi_k \circ \varphi_{i,k})^{-1}(v))) = \\
       \psi_k(\varphi_{i,k}(\varphi_{i,k}^{-1}(\varphi_k^{-1}(v)))) = \psi_k(\varphi_{k}^{-1}(v)).
   \end{gather*}
   Analogously we have that $\psi_i(\varphi_j^{-1}(v)) = \psi_k(\varphi_k^{-1}(v))$ and $\theta$ is well-defined. Consequently, $\varinjlim V_n^i = V_n$. Using Lemma \ref{lemma_entropy_limitedelcociente_cocientedellimite} and since $\varinjlim V_n^i= V_n$, we have that 
    \begin{equation}\label{rabiosa}
        \dim (V_n/V_{n-1}) = \sup_{i \in I}\left \{\dim \left (V_{n}^i/V_{n-1}^i \right ) \right \} 
    \end{equation}
    because the maps $\varphi_{i,j}$ are injective.
    
    For the last part, check that
    \begin{equation*}
    \frac{\log(\dim(V_n^i/V_{n-1}^i)}{n} \leq \limsup_{n\to \infty} \frac{\log(\dim(V_n^i/V_{n-1}^i)}{n}
    \end{equation*}
    for every $i \in I$ and every $n \geq 2$, then
    \begin{equation}
        \sup_{i\in I}\left \{  \frac{\log(\dim(V_n^i/V_{n-1}^i))}{n} \right \} \leq \sup_{i\in I} \left \{ \limsup_{n\to \infty} \frac{\log(\dim(V_n^i/V_{n-1}^i))}{n}\right \}
    \end{equation}
    for every $n \geq 2$. Finally, we have that
    
    \begin{gather*}
           \h(A,\F) = \limsup_{n\to \infty} \frac{\log \left ( \dim(V_n/V_{n-1})\right )}{n} = \\
           = \limsup_{n\to \infty} \frac{\log \left ( \sup_{i\in I} \left \{\dim(V_n^i/V_{n-1}^i) \right \}\right )}{n}  = \\
           = \limsup_{n \to \infty}  \sup_{i\in I}\left \{\frac{\log(\dim(V_n^i/V_{n-1}^i)}{n} \right \} \leq \sup_{i\in I} \left \{ \limsup_{n\to \infty} \frac{\log(\dim(V_n^i/V_{n-1}^i)}{n}\right \}= \\
           = \sup_{i\in I} \left \{\h(A_i,\F_i)\right \} = \h(A,S).
    \end{gather*}
\end{proof}

\section{Entropy of strict directed systems of finitely generated algebras}\label{sec_finitely_generated}

In this section, we address a consistency question 
concerning the equality of entropies computed through different standard direct systems of finitely filtered algebras whose underlying algebra limits agree. Then, we generalize some of the results in \cite{entropy2}.

\begin{definition}
We say that $S = (\{(A_i,\F_i)\}_{i \in I},\{\varphi_{i,j}\}_{j > i})$ is a \textit{standard direct system} if it is a strict direct system and every $K$-algebra $A_i$ is finitely generated. 
\end{definition}

\begin{proposition}\label{prop_entropy_independencia}
    Consider a $K$-algebra $A$ with $(A,\F) \cong \varinjlim (A_i,\F_i) \cong \varinjlim (B_{i'},\G_{i'})$ such that $S = (\{(A_i,\F_i)\}_{i \in I},\{\varphi_{i,j}\}_{j > i})$ and $S' = (\{(B_{i'},\G_{i'})\}_{i' \in I'},\{\psi_{i',j'}\}_{j' > i'})$ are two standard direct systems of finitely filtered algebras. Then for every $i \in I$ there is an $i' \in I'$ and a monomorphism $\omega_{i,i'} \colon A_i \to B_{i'}$ given by $\omega_{i,i'}(x) = z$ where $\varphi_i(x) = \psi_{i'}(z)$ which is strict. 
\end{proposition}
\begin{proof}
Without loss of generality, by the universal property we can consider that $A = \varinjlim A_i = \varinjlim B_{j'}$.
    Consider $i \in I$, since $A_i$ is finitely generated, take $\{a_1,\ldots, a_n\}$ generators, then $\varphi_i(A_i)$ is generated by $\{\varphi_i(a_1),\ldots, \varphi_i(a_n)\}$. Since $A = \varinjlim B_{i'}$ for each $\varphi_i(a_k)$ there is an $i_k' \in I'$ such that $\varphi_i(a_k) \in \psi_{i'_k}(B_{i'_k})$. The set $\{i'_1,\ldots, i'_n\}$ is finite and $I'$ is a directed set, then there is an upper bound $i' \in I'$ such that $i'_k \leq i'$ for every $k = 1,\ldots, n$. In particular, $\varphi_i(a_k) \in \psi_{i'}(B_{i'})$ for every $k=1,\ldots, n$. This is a direct consequence of $\varphi_i(a_k) = \psi_{i'_k}(z_k) = \psi_{i'}\psi_{i'_k,i'}(z_k)$. For each $a_k \in A_i$ we have that $b_k := \psi_{i'_k,i'}(z_k) \in B_{i'}$. Take $a \in A_i$, then $a = p(a_1,\ldots, a_n)$ for $p \in K\langle x_1,\cdots, x_n\rangle $, and
\begin{gather*}
       \varphi_i(a) = p(\varphi_i(a_1),\ldots,\varphi_i(a_k)) = p(\psi_{i'_1}(z_1),\psi_{i'_2}(z_2),\ldots, \psi_{i'_n}(z_n)) = \\  p(\psi_{i'}(b_1), \psi_{i'}(b_2),\ldots, \psi_{i'}(b_n)) = \psi_{i'}(p(b_1,\ldots,b_n))
\end{gather*} 
Summarizing, for every $a \in A_i$ there is a $b \in B_{i'}$ such that $\varphi_i(a) = \psi_{i'}(b)$. Let us see that such element is unique. If $\varphi_i(a) = \psi_{i'}(b') = \psi_{i'}(b)$, then $b=b'$ because $\psi_{i'}$ is injective. Next, we  can define

\begin{equation*}
    \begin{array}{cccc}
      \omega_{i,i'}\colon&  A_i & \to & B_{i'} \\
       & a & \mapsto & \omega_{i,i'}(a):= b 
    \end{array}
\end{equation*}
with $\varphi_i(a) = \psi_{i'}(b)$.

Let us prove that $\omega_{i,i'}$ is injective. Take $a, a'\in A_i$ with $\omega_{i,i'}(a) = \omega_{i,i'}(a')$ $\varphi_i(a) =  \psi_{i'}(\omega_{i,i'}(a)) = \psi_{i'}(\omega_{i,i'}(a'))= \varphi_i(a')$, by injectivity of $\varphi_{i}$ we have that $a = a'$.  

Next, let us prove that $\omega_{i,i'}$ is a linear map. Take $\alpha a_1 + a_2 \in A_i$, then there are $b_1,b_2 \in B_{i'}$ such that $\varphi_i(a_l) = \psi_{i'}(b_l)$ for $l =1,2$. Then $\varphi_i(\alpha a_1 +a_2) = \psi_{i'}(\alpha b_1 + b_2)$ which means that $\omega_{i,i'}(\alpha a_1+a_2) = \alpha b_1+b_2 = \alpha \omega_{i,i'}(a_1)+\omega_{i,i'}(a_2)$.

Analogously, we can prove that $\omega_{i,i'}(a_1 a_2) = \omega_{i,i'}(a_1)\omega_{i,i'}(a_2)$.

Now, we prove that $\omega_{i,i'}(V_n^i)\subset W_n^{i'}$ for
any $n$, $i$ and $i'$. If $a\in A_i\cap V_n^i$, then $\varphi_i(a)\in V_n=\bigcup_{i\in I} \varphi_i(V_n^i) / \sim $ but $\varphi_i(a)=\psi_{i'}(b)\in \psi_{i'}(B_{i'})\cap V_n$,
hence $b\in \psi_{i'}^{-1}(V_n)=\psi_{i'}^{-1}(W_n)=W_n^{i'}$.

Finally, we need to check that $\omega_{i,i'}^{-1}(W_n^{i'}) = V_n^i$. Consider $a \in \omega_{i,i'}^{-1}(W_n^{i'})$, then $b = \omega_{i,i'}(a) \in W_n^{i'}$. This element verifies that $\varphi_i(a) = \psi_{i'}(b) \in V_n$. Since $\varphi_i$ is strict by Lemma \ref{lemma_strict_morphism} we have that $a \in V_n^i$. 
\end{proof}

\begin{theorem}\label{th_independencia_entropía}
    Under the conditions of Proposition \ref{prop_entropy_independencia}, then 
    \begin{equation*}
        \h(A,S) = \h(A,S').
    \end{equation*}
\end{theorem}
\begin{proof}
    Consider $i \in I$, there is some $i' \in I'$ such that $\omega_{i,i'} \colon (A_i, \F_i) \to (B_{i'},\G_{i'})$ is a strict monomorphism. By  \cite[Lemma 4.2]{entropy1} we have that $\h(A_i, \F_i) \leq \h(B_{i'},\G_{i'}) \leq \h(A,S')$. By the properties of the supremum we have that $\h(A,S) \leq \h(A,S')$. Analogously we can prove that $\h(A,S') \leq \h(A,S)$.
\end{proof}

\begin{remark}\label{remark_entropia_sinsistemas}
     Consider $(A,\F)$ a filtered algebra such that $A$ is the direct limit of $(\{A_i\}_{i\in I},\{\varphi_{i,j}\}_{i \leq j})$, with $\varphi_{i,j}$ injective. In addition, if $A_i$ is finitely generated and $\varphi_i^{-1}(\F) = \F_i$ is a finite filtration of $A_i$, then $S = (\{(A_i,\F_i)\}_{i\in I}, \{\varphi_{i,j}\}_{i\leq j})$ is a standard direct system such that $(A,\F) = \varinjlim (A_i,\F_i)$. For algebras such that they verify these conditions, following Theorem \ref{th_independencia_entropía} it makes sense to construct a definition of entropy that does not depend on the system rather than the filtration that we consider for the whole algebra.
\end{remark}

The conditions given in Remark \ref{remark_entropia_sinsistemas} are verified by the path and Leavitt path algebras. 

\begin{definition}
    Consider $(A,\F)$ a filtered algebra under the conditions of Remark \ref{remark_entropia_sinsistemas}. Then we define:
    \begin{equation}
        \h(A,\F) = \sup_{i\in I}\{\h(A_i,\F_i)\}
    \end{equation}
    for a standard direct system of finitely filtered algebras. 
\end{definition}

\section{Entropy of path algebras and Leavitt path algebras}\label{sec_further_examples}
 In this section we specialise this notion of algebraic entropy to the context of path and Leavitt path algebras associated with infinite (row-finite) graphs.

\begin{remark}
Any infinite dimensional path algebra $KE$ is the direct limit of a strict finitely filtered system of finite dimensional path algebras. Indeed, we have $E=\varinjlim{S}$ where $S$ is a direct system in $\Dg$ of finite subgraphs of $E$. Then, by \cite[Lemma 2.5]{Goodearl} $KE=K(\varinjlim{S})=\varinjlim{KS}$ and $KS$ is a strict directed system of finitely filtered $K$-algebras (of finite dimension).
A similar result can be given for row-finite Leavitt path algebras \cite[Lemma 3.1 and 3.2]{AraMorenoPardo}.
\end{remark}

\begin{proposition}\label{prop_igualdad}
Let $E$ a row-finite infinite graph. The entropy of $KE$ and $L_K(E)$ coincides for the standard filtrations.
\end{proposition}
\begin{proof}
On the line of \cite[Lemma 3.2]{AraMorenoPardo} and \cite[Proposition 2.6]{Goodearl} the graph $E$ is the direct limit of its complete subgraphs $(\{E_i\}_{i\in I}, \{\mathfrak{i}_{i,j}\}_{i \leq j})$ where $i \leq j$ if $E_i \subseteq E_j$ is a complete subgraph, and $\mathfrak{i}_{i,j}$ is the inclusion of graphs. Furthermore, by \cite[Lemma 2.5]{Goodearl}, $L_K(E) = \varinjlim L_K(E_i)$ and $KE = \varinjlim KE_i$. Denote by $\F$, $\F_i$, $\G$ and $\G_i$ the standard filtrations for $KE$, $KE_i$, $L_K(E)$ and $L_K(E_i)$ respectively. Consider $S = (\{(KE_i,\mathcal F_i)\}_{i \in I},\{\varphi_{i,j}\}_{i \leq j})$ and $R = (\{(L_K(E_i),\mathcal G_i)\}_{i \in I},\{\psi_{i,j}\}_{i \leq j})$ the direct systems of filtered algebras induced by the inclusions of complete subgraphs. If $S$ and $R$ form two standard direct systems such that $(KE,\mathcal F) = \varinjlim (KE_i,\mathcal F_i)$ and $(L_K(E),\mathcal G) = \varinjlim (L_K(E_i),\mathcal G_i)$, then, by the Definition \ref{def_entropy_limit} and \cite[Theorem~3.5]{entropy2}, it follows that 
\begin{equation*} \h(KE,S) = \sup_{i\in I}\{\h(KE_i,\mathcal F_i)\} = \sup_{i\in I}\{\h(L_K(E_i),\mathcal G_i)\} = \h(L_K(E),R). \end{equation*} 
Therefore, it only remains to prove that $S$ and $R$ are standard direct systems such that $(KE,\mathcal F)=\varinjlim S$ and $(L_K(E),\mathcal G)=\varinjlim R$. The following diagram illustrates the situation for $i\leq j$:

\[
\begin{tikzcd}
KE_i\arrow[r,"\varphi_{i,j}"]\arrow[d,hook] & KE_j\arrow[r,"\varphi_j"]\arrow[d,hook] & KE\arrow[d,hook]\\
    L_K(E_i)\arrow[r,"\psi_{i,j}"'] & L_K(E_j)\arrow[r,"\psi_j"'] & L_K(E).\\
\end{tikzcd}
\]
By the definition of the standard filtration of path and Leavitt path algebras and being $\varphi_{i,j}$ and $\psi_{i,j}$ induced by the inclusion, it is trivial to check that both of them are strict monomorphisms of filtered algebras. 
Consider $\lambda \in V_n$, that is, $|\lambda | \leq n$, since $E$ is the direct limit of its complete subgraphs, there is some $i \in I$ such that $\lambda \in \PP(E_i)$ with $\lambda \in V_n^i \in \F_i$, which leads to $V_n = \cup_{i\in I} \varphi_i(V_n^i)$.
In addition, take $\lambda\mu^* \in W_n$. Since $E$ is the directed limit of its complete subgraphs, there is $k,j \in I$ such that $\lambda \in \PP(E_j)$ and $\mu \in \PP(E_k)$, since the complete subgraphs form a directed set, there is an index $i \in I$ with $i \geq k,j$ such that $\lambda,\mu \in \PP(E_i)$, then $\lambda \mu^* \in W_n^i$. Finally this implies that $W_n = \bigcup_{i \in I} \psi_{i}(W_n^i)$. Consequently, we have that $(KE,\F) = \varinjlim (KE_i,\F_i)$ and $(L_K(E),\G) = \varinjlim (L_K(E_i), \G_i)$.
\end{proof}

To end this section, we consider several examples which show that the  
calculus may have a twisted prong when dealing with  infinite dimensional path algebras. In particular, there is no systematic procedure to compute the entropy since the calculus leads one to many different scenarios.

\begin{example}
Consider the directed graph $E$ consider in Figure \ref{fig_morcilla_infinita}.
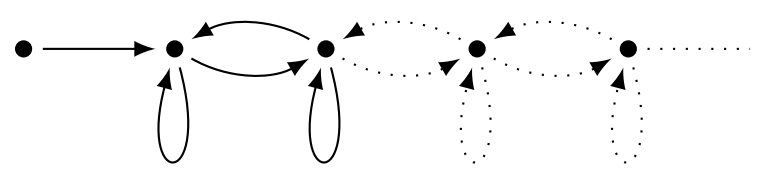
\begin{figure}[ht]
    \centering
   \begin{tikzpicture}
[scale=0.5,shorten <=2pt,shorten >=2pt,>=latex, node distance={15mm}, main/.style = {draw, fill, circle, inner sep = 1pt}]

\node[label=above:$E:$] (-1) at (-2,-.5) {};

\node[main] (0) at (0,0) {};
\node[main] (1) at (2,0) {};
\node[main] (2) at (4,0) {};
\node[main] (3) at (6,0) {};
\node[main] (4) at (8,0) {};
\node (5) at (10,0) {};

\draw[->] (1) edge[loop below, looseness = 80] (1);
 \draw[->] (2) edge[loop below, looseness = 80] (2);
\draw[->,dotted] (3) edge[loop below, looseness = 80] (3);
\draw[->,dotted] (4) edge[loop below, looseness = 80] (4);
  
\draw[->] (0) to (1);
\draw[->] (1) [bend right] to (2);
\draw[->] (2) [bend right] to (1); 
\draw[->,dotted] (2) [bend right] to (3);
\draw[->,dotted] (3) [bend right] to (2);
\draw[->,dotted] (4) [bend right] to (3);
\draw[->,dotted] (3) [bend right] to (4);
\draw[dotted] (4) to (5);
\end{tikzpicture}
\caption{Directed graph $E$.}
\label{fig_morcilla_infinita}
\end{figure}

With the aim of computing the entropy of $KE$, we define the finite subgraphs $E_n$ displayed in Figure \ref{doscosas}, together with their adjacency matrices 

      \[
        A_n= \tiny
        \begin{pmatrix}
        0 & 1 & 0      & 0      & \cdots & 0 &0\\
        0 & 1 & 1      & 0      &\cdots   &0  &  0\\
        0 & 1 & 1      & 1      & \cdots  &0 & 0\\
        0 & 0 & 1 & 1 & \cdots   & 0& 0 \\
         \vdots  & \vdots   & \vdots & \vdots & \ddots       & \vdots  & \vdots\\
         0 & 0 & 0 & 0 & \cdots  & 1 & 1 \\
        0 & 0 & 0  & 0      & \cdots       & 1  & 1  
        \end{pmatrix}.
        \]

\begin{figure}[ht]
    \centering
        \centering
        \begin{tikzpicture}[
            scale=0.5,
            shorten <=2pt,
            shorten >=2pt,
            >=latex,
            main/.style={draw, fill, circle, inner sep=1pt}
        ]
        \node[main, label=above:$v_0$] (0) at (0,0) {};
        \node[main, label=above:$v_1$] (1) at (2,0) {};
        \node[main] (2) at (4,0) {};
        \node[main] (3) at (6,0) {};
        \node[main, label=above:$v_n$] (4) at (8,0) {};
        \draw[->] (1) edge[loop below, looseness=80] (1);
        \draw[->] (2) edge[loop below, looseness=80] (2);
        \draw[->] (3) edge[loop below, looseness=80] (3);
        \draw[->] (4) edge[loop below, looseness=80] (4);
        \draw[->] (0) to (1);
        \draw[->] (1) [bend right] to (2);
        \draw[->] (2) [bend right] to (1);
        \draw[->,dotted] (2) [bend right] to (3);
        \draw[->,dotted] (3) [bend right] to (2);
        \draw[->] (4) [bend right] to (3);
        \draw[->] (3) [bend right] to (4);
\end{tikzpicture}
        \caption{Subgraph $E_n$ with adjacency matrix $A_n$.}
        \label{doscosas}
\end{figure}
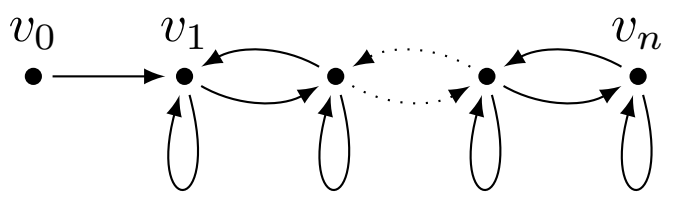

The algebraic entropy of $KE_n$ relative to its standard filtration, denoted by $\F_n$, is $\log\left (1+2\cos\frac{\pi}{n+1}\right)$. Indeed, applying \cite[Corollary 5.5]{entropy1} the algebraic
 entropy of $KE_n$ is the same as the one of $K(E_n\setminus\{{v_0}\})$. So we can remove the source $v_0$
 and we have the graph  
 
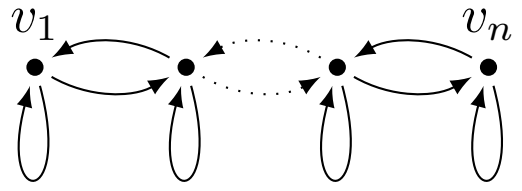
\[
        \begin{tikzpicture}[
            scale=0.5,
            shorten <=2pt,
            shorten >=2pt,
            >=latex,
            main/.style={draw, fill, circle, inner sep=1pt}
        ]
%
        \node[main, label=above:$v_1$] (1) at (2,0) {};
        \node[main] (2) at (4,0) {};
        \node[main] (3) at (6,0) {};
        \node[main, label=above:$v_n$] (4) at (8,0) {};
        \draw[->] (1) edge[loop below, looseness=80] (1);
        \draw[->] (2) edge[loop below, looseness=80] (2);
        \draw[->] (3) edge[loop below, looseness=80] (3);
        \draw[->] (4) edge[loop below, looseness=80] (4);
%
        \draw[->] (1) [bend right] to (2);
        \draw[->] (2) [bend right] to (1);
        \draw[->,dotted] (2) [bend right] to (3);
        \draw[->,dotted] (3) [bend right] to (2);
        \draw[->] (4) [bend right] to (3);
        \draw[->] (3) [bend right] to (4);
\end{tikzpicture}
\]
The adjacency matrix of this graph is a tridiagonal matrix $n\times n$, that is a matrix whose nonzero entries are $1$ in the three diagonals. In \cite[Section 26.5, Tridiagonal Toeplitz Matrices, p. 524]{higham} it is proven that the eigenvalues of this matrix are $1+2\cos\frac{k\pi}{n+1}$ with $k = 1, \ldots, n$. Thus, by \cite[Theorem 5.4]{entropy1}
\begin{equation*}
    \h(KE_n,\F_n)=\log \left ( \max_{k = 1, \ldots,n} \left \{1+2\cos\frac{k\pi}{n+1} \right \}\right ).
\end{equation*}

If we consider the standard system $S$ induced by the inclusion of path algebras $\{(KE_n,\F_n)\}_{n\geq 1}$.
Then applying Definition \ref{def_entropy_limit}, we have that

\begin{gather*}
        \h(KE,S) = \sup_{n\geq 1}\left \{\log \left ( \max_{k = 1, \ldots,n} \left \{1+2\cos\frac{k\pi}{n+1} \right \}\right ) \right \} = \\ \sup_{n\geq 1} \left \{\log\left (1+2\cos\frac{\pi}{n+1}\right )\right \}=\log 3.
\end{gather*}

\end{example}

\medskip

\noindent\textbf{Data Availability Statement:} The authors confirm that the data supporting the findings of this study are available within the article.

\medskip

\noindent\textbf{Authorship contribution statement:} All authors contributed equally to this work.

\medskip

\noindent\textbf{Ethical approval:} This declaration is not applicable.

\medskip

\noindent\textbf{Declaration of competing interest:} The authors declare that they have no competing interests.

\noindent\textbf{Funding:} The authors are supported by the Spanish Ministerio de Ciencia e Innovaci\'on y Universidades and by the Junta de Andaluc\'{i}a through projects PID2023-152673NB-I00 and PPRO-FQM336-G-2023 (FQM336-G-FEDER) respec- tively.

\bibliographystyle{plain} 
\bibliography{ref.bib}

@article{entropy1,
  title={Algebraic Entropy of Path Algebras and Leavitt Path Algebras of Finite Graphs.},
  author={Bock, Wolfgang and Canto, Crist{\'o}bal Gil and Martín Barquero, Dolores and Martín Gonz{\'a}lez, C{\'a}ndido and Ruiz Campos, Iv{\'a}n 
  and Sebandal, Alfilgen},
  journal={Results in Mathematics},
  volume={79},
  number={5},
  pages={180},
  year={2024},
  publisher={Springer}
}

@article{entropy2,
  title={The Algebraic Entropies of the Leavitt Path Algebra and the Graph Algebra Agree.},
  author={Bock, Wolfgang and Canto, Crist{\'o}bal Gil and Martín Barquero, Dolores and Martín Gonz{\'a}lez, C{\'a}ndido and Ruiz Campos, Iv{\'a}n 
  and Sebandal, Alfilgen},
  journal={Results in Mathematics},
  volume={79},
  number={8},
  pages={266},
  year={2024},
  publisher={Springer}
}

@book{maclane,
  author    = {Mac Lane, Saunders},
  title     = {Categories for the Working Mathematician},
  series    = {Graduate Texts in Mathematics},
  volume    = {5},
  edition   = {Second},
  year      = {1998},
  publisher = {Springer-Verlag},
  address   = {New York},
  doi       = {10.1007/978-1-4757-4721-8},
  isbn      = {978-0-387-98403-2}
}

@incollection{Goodearl,
  author    = {Goodearl, Kenneth R.},
  title     = {Leavitt path algebras and direct limits},
  booktitle = {Rings, Modules, and Representations},
  editor    = {Dung, Nguyen Viet and Facchini, Alberto and Wisbauer, Robert and Zhou, Yiqiang},
  series    = {Contemporary Mathematics},
  volume    = {480},
  pages     = {165--188},
  year      = {2009},
  publisher = {American Mathematical Society},
  address   = {Providence, Rhode Island},
  doi       = {10.1090/conm/480/09440}
}

@book{higham,
  author    = {Nicholas J. Higham},
  title     = {Accuracy and Stability of Numerical Algorithms},
  edition   = {Second},
  publisher = {Society for Industrial and Applied Mathematics (SIAM)},
  address   = {Philadelphia, PA, USA},
  year      = {2002},
  isbn      = {0-89871-521-0}
}

@article{AbramsPino, title={The \uppercase{L}eavitt path algebra of a graph}, author={Abrams, Gene  and Aranda Pino, Gonzalo}, journal={Journal of  Algebra}, pages={319--334}, volume={293}, year={2005}}

@article{AraMorenoPardo,
author={Ara, Pere and Moreno, Mar{\'i}a A. and Pardo, Enrique}, title={Nonstable {K-theory} for graph algebras}, year={2007}, pages={157--178},
journal={Algebras and Representation Theory},
volume={10}
}

@article{dikranjan2012connection,
  title={The connection between topological and algebraic entropy},
  author={Dikranjan, Dikran and Bruno, Anna Giordano},
  journal={Topology and its Applications},
  volume={159},
  number={13},
  pages={2980--2989},
  year={2012},
  publisher={Elsevier}
}

@article{dikranjan2016entropy,
  title={Entropy on abelian groups},
  author={Dikranjan, Dikran and Bruno, Anna Giordano},
  journal={Advances in Mathematics},
  volume={298},
  pages={612--653},
  year={2016},
  publisher={Elsevier}
}

@article{bruno2017some,
  title={Some properties of the growth and of the algebraic entropy of group endomorphisms},
  author={Bruno, Anna Giordano and Spiga, Pablo},
  journal={Journal of Group Theory},
  volume={20},
  number={4},
  pages={763--774},
  year={2017},
  publisher={De Gruyter}
}

@article{newman2000entropy,
  title={The entropy of graded algebras},
  author={Newman, Michael Frederick and Schneider, Csaba and Shalev, Aner},
  journal={Journal of Algebra},
  volume={223},
  number={1},
  pages={85--100},
  year={2000},
  publisher={Elsevier}
}

@book{Boltzmann,
  title={{\"U}ber die mechanische Bedeutung des zweiten Hauptsatzes der W{\"a}rmetheorie},
  author={Boltzmann, Ludwig},
  year={1866},
  publisher={Staatsdruckerei}
}

@article{Zelmanov12,
author = {Alahmadi, Adel and Alsulami, Hamed and Jain, Surender K. and Zelmanov, Efim},
title = {Leavitt path algebras of finite {Gelfand}-{Kirillov} dimension},
journal = {Journal of Algebra and Its Applications},
volume = {11},
number = {06},
pages = {1250225},
year = {2012},
doi = {10.1142/S0219498812502258},

URL = { https://doi.org/10.1142/S0219498812502258
},
eprint = { 
    
        https://doi.org/10.1142/S0219498812502258
    
}
,

}

@article{moreno2018graph,
  title={Graph algebras and the {Gelfand}--{Kirillov dimension}},
  author={Moreno-Fern{\'a}ndez, Jos{\'e} M. and Molina, Mercedes Siles},
  journal={Journal of Algebra and Its Applications},
  volume={17},
  number={05},
  pages={1850095},
  year={2018},
  publisher={World Scientific}
}

@misc{Mathematica,
  author = {Wolfram Research{,} Inc.},
  title = {Mathematica, {V}ersion 13.1},
  url = {https://www.wolfram.com/mathematica},
  note = {Champaign, IL, 2022}
}

@article{Adler,
  title={Topological entropy},
  author={Adler, Roy L. and Konheim, Alan G. and McAndrew, M. Harry},
  journal={Transactions of the American Mathematical Society},
  volume={114},
  number={2},
  pages={309--319},
  year={1965},
  publisher={JSTOR}
}

@article{Weiss,
  title={Algebraic and other entropies of group endomorphisms},
  author={Weiss, Michael D.},
  journal={Mathematical systems theory},
  volume={8},
  number={3},
  pages={243--248},
  year={1974},
  publisher={Springer}
}

@article{peters1979,
  title={Entropy on discrete abelian groups},
  author={Peters, Justin},
  journal={Advances in Mathematics},
  volume={33},
  number={1},
  pages={1--13},
  year={1979},
  publisher={Academic Press}
}

@misc{stanley,
  title={Enumerative combinatorics. Volume 2. Paperback ed},
  author={Stanley, Richard P},
  year={2001},
  publisher={Cambridge: Cambridge University Press, paperback ed. ed}
}

@book{weibel1994introduction, 
place={Cambridge},
series={Cambridge Studies in Advanced Mathematics}, title={An Introduction to Homological Algebra}, publisher={Cambridge University Press}, 
author={Weibel, Charles A.}, 
year={1994}, 
collection={Cambridge Studies in Advanced Mathematics}}

@article{entropydaniel,
    author = {Gonçalves, Daniel and Royer, Danilo and Augusto Tasca, Felipe},
    title = "{Entropy of Local Homeomorphisms With Applications to Infinite Alphabet Shift Spaces}",
    journal = {International Mathematics Research Notices},
    volume = {2024},
    number = {6},
    pages = {4842-4892},
    year = {2023},
    month = {08},
    issn = {1073-7928},
    doi = {10.1093/imrn/rnad205},
    url = {https://doi.org/10.1093/imrn/rnad205},
    eprint = {https://academic.oup.com/imrn/article-pdf/2024/6/4842/56930491/rnad205.pdf},
}

@article{entropychinos,
  title={Entropy dimension of shift spaces on monoids},
  author={Ban, Jung-Chao and Chang, Chih-Hung and Huang, Nai-Zhu},
  journal={Journal of Mathematical Physics},
  volume={61},
  number={7},
  year={2020},
  publisher={AIP Publishing}
}
\end{document}